\documentclass[12pt,oneside,reqno]{amsart}

\usepackage[T1]{fontenc}
\usepackage[utf8]{inputenc}
\usepackage{amsmath}
\usepackage{lipsum}
\usepackage{amssymb}
\usepackage{amsfonts}
\usepackage{enumitem}
\usepackage{stackengine}
\usepackage{graphicx}
\usepackage{calc}

\usepackage{hyperref}
\usepackage{ragged2e}
\usepackage{bbm}
\usepackage{xcolor}
\hypersetup{colorlinks=true}
\usepackage[a4paper, left=2.5cm, right=2.5cm, top=2.5cm, bottom=2.5cm]{geometry}
\usepackage{amsthm}
\usepackage{indentfirst}
\makeatletter
\DeclareMathAlphabet{\mathcal}{OMS}{zplm}{m}{n}

\newcommand{\s}{{\small{$\Sigma$}}}

\newcommand{\tmod}{\mathrm{mod}}
\newcommand{\x}{\mathbf{x}}

\newcommand{\y}{\mathbf{y}}

\newcommand{\R}{\mathbb{R}}
\newcommand{\T}{\mathbb{T}}

\newcommand{\Q}{\mathbb{Q}}
\newcommand{\N}{\mathbb{N}}

\newcommand{\E}{\mathbb{E}}

\newcommand{\Z}{\mathbb{Z}}

\newcommand{\cH}{\mathcal{H}}
\newcommand{\cL}{\mathcal{L}}

\newcommand{\cA}{\mathcal{A}}

\newcommand{\dimH}{\dim_{\mathrm{H}}}
\newcommand{\dist}{\mathrm{dist}}

\newcommand{\linte}[1]{{\left\lfloor #1
		\right\rfloor}}

\newcommand{\one}{\mathbbm{1}}

\newcommand{\PP}{\mathcal P}
\newcommand{\all}{\;\;\forall}

\theoremstyle{plain}

\newtheorem{theorem}{Theorem}
\newtheorem{lemma}{Lemma}

\newtheorem{corollary}{Corollary}

\newtheorem*{claim*}{Claim}

\theoremstyle{definition}

\newtheorem{example}{Example}

\newtheorem*{examples*}{Examples}
\newtheorem*{example*}{Example}

\newtheorem*{notations*}{Notations}
\newtheorem*{notation*}{Notation}

\theoremstyle{remark}
\newtheorem{remark}{{Remark}}

\newtheorem{lemmaCBC}{{\rm \textbf{Lemma CBC (Convergence Borel--Cantelli)}} \!\!\!\!}

\newtheorem{lemmaDBC}{{\rm \textbf{Lemma DBC (Divergence Borel--Cantelli)}} \!\!\!\!}

\newtheorem{lemmaQBC}{{\rm \textbf{Lemma QBC (Quantitative Borel--Cantelli)}} \!\!\!\!}

\numberwithin{equation}{section}
\numberwithin{question}{section}
\numberwithin{theorem}{section}
\numberwithin{thm}{section}
\numberwithin{lemma}{section}
\numberwithin{proposition}{section}
\numberwithin{cor}{section}
\numberwithin{corollary}{section}
\numberwithin{claim}{section}
\numberwithin{definition}{section}
\numberwithin{example}{section}
\numberwithin{remark}{section}
\numberwithin{notations}{section}
\numberwithin{notation}{section}
\numberwithin{claim}{section}
\numberwithin{problem}{section}
\numberwithin{figure}{section}

\allowdisplaybreaks
\renewcommand{\emph}[1]{{\it#1}}

\newtheorem{claimf}{{\rm \textbf{Claim~F}}  \!\!\!\!\!\!}

\begin{document}
\title{Shrinking Targets versus Recurrence: a brief survey}

\author[Y. He]{Yubin He}
\address[Y. He]{Department of Mathematics, Shantou University, Shantou, Guangdong, 515063,
China}
\email{ybhe@stu.edu.cn}

\author[B. Li]{Bing Li}
\address[B. Li]{School of Mathematics, South China University of Technology,  Wushan Road 381, Tianhe District, Guangzhou, China}
\email{scbingli@scut.edu.cn}

\author[S. Velani]{Sanju Velani}
\address[S. Velani]{Department of Mathematics,
University of York, Heslington, York, YO10
5DD, England.}
\email{sanju.velani@york.ac.uk}


\begin{abstract}

Let $(X,d)$  be a compact metric space and $(X,\mathcal{A},\mu,T)$ a measure preserving dynamical system.  Furthermore, given a real, positive function
$\psi$,  let $W(T, \psi)$ and
$
R(T,\psi) $  respectively
denote the shrinking target set  and the recurrent set associated with the dynamical system.   Under certain mixing properties it is known that if  the natural measure sum 
diverges then the recurrent and shrinking target sets are of full $\mu$-measure.  The purpose of this survey is to provide  a brief  overview of such results, to discuss the potential quantitative strengthening of the full measure statements and to bring to the  forefront  key differences in the theory.
\end{abstract}
\maketitle

\section{Introduction: background and motivation}  \label{bgm}

Let $(X,d)$  be a compact metric space and $(X,\mathcal{A},\mu,T)$ be an ergodic probability measure preserving system.  Furthermore, given a real, positive function $\psi:\N\to\R_{\ge 0}$ let
$$
R(T,\psi) := \big\{ x \in X :  T^nx \in B(x, \psi(n)) \hbox{ for infinitely many }n\in \N  \big\}  $$
denote the associated \emph{recurrent set}, and given a point  $x_0 \in X$ let
$$
W( T, \psi) := \big\{ x \in X : T^nx \in B(x_0, \psi(n))  \hbox{ for infinitely many }n\in \N  \big\}  $$
denote the  associated \emph{shrinking target set}.   If $\psi=c$ (a constant),   it follows from  two  foundational
results in dynamics, namely the Poincar\'{e} Recurrence  Theorem and  the Ergodic Theorem (see \cite[Theorem 1.4 and Theorem 1.14]{walters2000}),  that
$$
  \mu( R(T,c) )   = 1 =
 \mu( W(T,c) )  \,    .  $$
 \noindent Note that we do not need  the system to be  ergodic to conclude that  $\mu( R(T,c) )   = 1$.  The upshot is that in both the  shrinking target and recurrence setups the trajectories of almost all points
 will hit  `constant' balls  infinitely often.   In view of this, it is natural to ask:
what is the  $\mu$-measure of  the sets  if $\psi(n) \to 0 $ as $n \to \infty$?
In turn, whenever the $\mu$-measure is zero,  it is natural to ask about the Hausdorff dimension of the sets under consideration. Both these questions fall under the general umbrella of the ``shrinking target problem'' formulated in \cite{hill1995}.   In this survey we will concentrate mainly on the measure question.

To start with let us deal with  problem of determining ``natural''  conditions under which the $\mu$-measure of the sets under consideration is zero.  This turns out to be relatively straightforward once we observe that both  $  R( T, \psi)   $    and $  W( T, \psi) $ are $\limsup$ sets.  Indeed, by definition
$$
R( T, \psi) =  \limsup_{n \to \infty} R_n
$$
where for each $n\in \N $
\begin{eqnarray}\label{def_An}
 R_n \,  =   \, R_n( T,\psi)  & := & \big\{ x \in X :  T^nx \in B(x, \psi(n)) \big\}
  \, .
\end{eqnarray}
Similarly, by definition
$$
W( T, \psi) =  \limsup_{n \to \infty} W_n
$$
where for each $n\in \N $
\begin{eqnarray}
 W_n \,  =   \, W_n( T,\psi)
  & := & \big\{ x \in X :  T^nx \in B_n:=B(x_0, \psi(n)) \big\}  \nonumber \\[1ex]  & = &  T^{-n} (B_n) \label{sousefulintro} \, .
\end{eqnarray}
From the onset, it is worth highlighting the obvious fact that  for   the shrinking target set $W( T, \psi)$ the ``targets''   $B_n$ do not depend on the  initial  point $x$  of the orbit $\{T^nx\}_{n \in \N}$.  As we shall see later on, this and in turn the fact that $W_n$ can be written as the pre-image of the target ball $B_n$ makes life significantly easier when considering the shrinking target setup.
For the moment, we simply observe that in both setups, a straightforward consequence  of the  $\limsup$ nature of the sets and the  (convergence) Borel–Cantelli Lemma (see \S\ref{theBCsec}) that:
\begin{equation}   \label{rques}
          \bullet  \ \hspace{2ex} \ \mu(R(T,\psi)) = 0  \quad \ {\rm if \quad }  \textstyle {\sum_{n=1}^{\infty}} \; \mu(R_n) \; <\infty \ ,
\end{equation}
and
  \begin{equation}  \label{stques}
          \bullet \hspace{2ex} \ \mu(W(T,\psi)) = 0  \quad \ {\rm if \quad }  \textstyle {\sum_{n=1}^{\infty}} \; \mu(W_n) \; <\infty \ .
\end{equation}\\
Thus, in both setups we derive a simple zero-measure criterion based on the convergence of the natural measure sum, which involves the ``building blocks'' of the $\limsup $  set. Now concentrating our attention to the measure (rather than dimension) aspect of the shrinking target problem,  the obvious question that arises at this point is: \!\!\!\!  \emph{what happens when the measure sum diverges?}
 In attempting to answer this question, we uncover some interesting new phenomena (both in terms of measure and counting) that appear within the recurrence setup but not in the shrinking target setup. In fact, in the latter setup, the results obtained when addressing this question align with findings from related areas — not only in ergodic theory and dynamical systems, but also in probability theory  \cite{BeresnevichVelani7,bill} and metric number theory \cite{BRV2016,Harman1998,khintchine1924}.   From this perspective, the $\limsup$ sets associated with the  recurrence setup prove to be far more intriguing.   The main body of work related to the recurrence setup has emerged relatively recently, with most developments occurring over the past five years \cite{allen2025,baker2021,baker2024,chang2019,he2024,he2022,jungie,hussain2022,Kirsebom2023,kleinbock2023,persson2023,persson2025,Sponheimer2025}. The discoveries leading to the new phenomena mentioned above are even more recent, having taken shape within the past year or so. The aim of this survey is to provide a basic overview of these developments and to highlight  some obvious open problems that remain both intriguing and ripe for further exploration.

The structure of this survey is as follows.  In   \S\ref{theBCsec} we  briefly discuss various forms of the classical  Borel--Cantelli Lemma in probability theory.  These statements for general $\limsup$ sets  provide a guide to the type of measure and counting results we could expect to prove for the recurrent and shrinking target $\limsup$ sets.       In \S\ref{STSU} we show that  the ``expected'' results hold for the shrinking target set. More precisely, we show that the results are essentially a direct consequence of exponentially mixing -- a basic assumption on the dynamical system  that  underpins the developments  within the recent   recurrence setup.  The latter  is the subject  of   \S\ref{RSU} and  things get rather interesting!   In short, the recurrence theory is  significantly  richer.

\vspace*{4ex}

\noindent{\emph{Acknowledgements.}}   SV would like sincerely to thank the organisers of the  conference `Fractal Geometry and Stochastics 7' for  the opportunity to give a talk  --  it was a stimulating and thoroughly  enjoyable experience.  The subject matter of that talk and my interactions with various participants at the conference  has formed the foundations for this work.  On a different note, I would also like to thank Peter ``Der Kaiser'' Selby -- a legendary figure widely regarded as one of the greatest footballers ever to grace the Theatre of Dreams at Fearnville.   His energy and desire to ``make it happen'' week in week out are both impressive and much appreciated. Finally, I’m touched that Ayesha and Iona still come to play from time to time — and even pass to this old git! BL was supported by National Key R\&D Program of China (No. 2024YFA1013700), NSFC 12271176 and Guangdong Natural Science Foundation 2024A1515010946. YH was supported by NSFC 12401108.  We would to thank the referee for their extremely detailed and useful comments which  have greatly improved the  survey.

\section{Zero-One measure criteria  for general lim-sup sets} \label{theBCsec}

 To set the scene, let $(\Omega,\cA,\mu)$ be a probability space  and
let $A_n$ ($n \in \mathbb{N}$) be a family of  measurable subsets (events) of $\Omega$.  Also, let
$$
A_{\infty}:=\limsup_{n \to \infty} A_n := \bigcap_{t=1}^{\infty}
\bigcup_{n=t}^{\infty} A_n \ ;
$$
i.e. $ A_{\infty} $ is the set of $x \in \Omega$ such that
$ x \in A_n $ for infinitely many $n \in \N$.

Determining the measure  of $A_{\infty}$ turns out to be one of the fundamental problems considered within the framework of classical probability theory -- see for example   \cite[Chp.1 \S4]{bill}   and    \cite[Chp.47]{Port}                                                                                                                      for general background and further details. With this in mind, the following \emph{convergence  Borel--Cantelli Lemma}\/ provides a beautiful and truly simple criterion for zero measure.

\bigskip

\begin{lemmaCBC}
 \emph{Suppose that
  $ \sum_{n=1}^{\infty} \mu(A_n) <~\infty$.  Then,
  $$\mu(A_{\infty}) = 0 \, . $$}
\end{lemmaCBC}

\medskip

\noindent  This powerful lemma, which is also known as the \emph{first  Borel--Cantelli Lemma},  has applications in numerous disciplines. In particular, within the context of number theory it is very much at the heart of Borel's proof that almost all numbers are normal \cite{borel}.

In view of  Lemma CBC, it is natural to ask whether or not there is a sufficient condition that enables us to deduce that the measure of $A_{\infty}$ is positive or possibly even full; that is to say that
$$\mu(A_{\infty}) = \mu(\Omega)=1   \, . $$

\noindent  The divergence of the measure sum $\sum_{n=1}^{\infty} \mu(A_n)$ is clearly necessary but certainly not enough as the following simple example demonstrates.

\medskip

\noindent \emph{Example}.  For $n \in \N$, let  $A_n=(0, \frac1n) \subset \Omega := [0,1]$ and $ \mu$ be one-dimensional Lebesgue measure restricted to $[0,1]$. Then
$$
\textstyle{\sum_{n=1}^\infty \mu(A_n)=\sum_{n=1}^\infty n^{-1} =\infty } $$
but
$$ A_{\infty}= \textstyle{ \bigcap_{t=1}^\infty \bigcup_{n=t}^\infty } A_n  =\bigcap_{t=1}^\infty (0,\tfrac1t)=\varnothing  \qquad {\rm and \ so  \ } \quad  \mu(A_{\infty}) = 0 \, .
$$

\medskip

The problem in the above example is that the building blocks $A_n$ of the $\limsup$ set under consideration  overlap `too much' - in fact they are nested.     The upshot is that in order to have  $\mu(A_{\infty}) > 0$,   we not only need the sum  of the measures to diverge but also  that the sets $A_n$ are in `some sense' independent; that is, we need to control overlaps!
Indeed, in early part of the last century, Borel $\&$ Cantelli essentially showed that  \textit{pairwise independence}   in the classical probabilistic  sense, which means that
\begin{equation}\label{fullind}
\mu ( A_s \cap A_t  )  =  \mu(A_s) \mu(A_t)   \qquad    s \neq t  \, , 
\end{equation}
implies that $\mu( A_{\infty}) = 1 $.  We say ``essentially'' since they actually required the building blocks $A_n$  to be   pairwise independent. In any case, not much later, it was shown that under the hypothesis of their  full measure statement, often referred to as the \emph{second Borel--Cantelli Lemma}, it is  possible to obtain a  significantly stronger quantitative statement.   In short,  given $x \in \Omega$  and  $ N \in \N$,    consider the \textit{counting function} $A(x,N)$ that counts the number of integers $n \le N$ such that $ x \in A_n$; that is
\begin{eqnarray} \label{countdef}
  A(x,N)  \, :=  \,    \# \big\{ 1\le n \le   N :    x \in A_n  \big\}\, .
\end{eqnarray}
Then,  if $ \sum_{n=1}^{\infty} \mu(A_n) = ~\infty$  and \eqref{fullind} holds, we have that
\begin{equation}  \label{asytriv}
A(x,N)  \ \sim  \   \sum_{n=1}^N \mu(A_n)     \qquad {\rm as  }  \qquad  N \to \infty \, ,
\end{equation}
for $\mu$-almost all $x\in \Omega$\footnote{{\small 
Throughout, given  functions $f$ and $g$  defined on a set $S$, we write $f(x) \sim g(x) $ as $ x \to x_0$ if  $\lim_{x \to x_0} f(x)/g(x) = 1 $.
Also, we 
write $f\ll g$ if there exists a constant $\kappa=\kappa(f,g,S)>0$, such that $|f(x)|\leq\kappa |g(x)|$ for all $x\in S$, and we write $f\asymp g$ if $f\ll g\ll f$.}}.  In other words, pairwise independence is in fact a  strong enough condition to describe  the asymptotic behaviour of the counting function  whereas the  second Borel--Cantelli Lemma simply tells us that
$$
 A(x,N)\;  \to  \;  \infty \qquad {\rm as  }  \qquad  N \to \infty \, ,
$$  for $\mu$-almost all $x\in \Omega$.    We will see in \S\ref{mechCR}  that  the asymptotic statement  \eqref{asytriv} can be strengthened.

In view of the above discussion, it seems that we are in good shape when the measure sum $\sum_{n=1}^{\infty} \mu(A_n)$    diverges   - we have a sufficient condition (namely  \eqref{fullind}) that  not only  guarantees full measure  but a  significantly stronger quantitative statement.      However, there is a serious downside to the  second Borel--Cantelli Lemma.  In many applications, we rarely have pairwise independence.   What is much more useful is the following variant which these days is often referred to as the \emph{divergence Borel--Cantelli Lemma}.

\begin{lemmaDBC}
 \!
\emph{Suppose that $ \sum_{n=1}^{\infty} \mu(A_n) =\infty$
and that there exists a constant $C\ge 1 $ such that
\begin{equation}\label{vbx1x}
\sum_{s,t=1}^Q  \mu(A_s\cap A_t)\le C\left(\sum_{s=1}^Q  \mu(A_s)\right)^2
\end{equation}
holds for infinitely many $Q\in\N$.
Then
$$
\mu(A_{\infty}) \ge C^{-1}\,.
$$
In particular, if $C=1$ then $\mu(A_{\infty}) =1 \, . $}
\end{lemmaDBC}

\medskip

We refer the reader to  \cite{Harman1998, Port, Sprindzuk} for the proof of the lemma which is essentially a consequence of the Cauchy--Schwarz inequality.  Condition \eqref{vbx1x} is often refereed to as \emph{quasi-independence on average} and  together with the divergence of the measure sum guarantees that the associated  $\limsup$  set $A_{\infty}$ is of positive measure. It does not in general guarantee full measure. However, this is not an issue if we already know by some other means (such as Kolmogorov's theorem \cite[Theorems 4.5 \& 22.3]{bill} or ergodicity \cite[\S24]{bill}) that the $\limsup $ set  $A_{\infty}$ satisfies a \emph{zero-one law}; namely  that
\begin{equation*}
\mu(A_{\infty}) =   \   0 \quad\text{or} \quad 1.
\end{equation*}
Alternatively, without the presence of a general zero-one law,  if we are willing to impose a little more structure on the probability space and \eqref{vbx1x} holds  `locally',   we can guarantee full measure.    More precisely,   assuming $\Omega$ is equipped with a metric such that
$\mu$ becomes a doubling Borel measure, we can guarantee that
$\mu(A_{\infty}) = \mu(\Omega)=1$ if we can establish \emph{local quasi-independence on average}\index{local quasi-independence on average}; that is, we replace
\eqref{vbx1x} in Lemma~DBC by the condition that
\begin{equation}\label{vbx1xslv}
  \sum_{s,t=1}^Q  \mu\big((B \cap A_s) \cap (B \cap A_t) \big)\le \frac{C}{\mu(B)} \left(\sum_{s=1}^Q  \mu(B \cap A_s)\right)^2 \, \smallskip
\end{equation}
for any sufficiently small ball $B:=B(x, r)$ with center $x$ in $\Omega$ and
$\mu(B) > 0$. The constant $C \ge 1 $ is independent of the ball $B$.
Recall, that $\mu $ is said to be \emph{doubling}
if there are constants $ \lambda  \geq 1$ and $r_0>0$ such that for any $x \in \Omega$  and $0<r<r_0$
\begin{equation}\label{doub}
\mu (B(x, 2r))  \, \leq \, \lambda \ \mu(B(x,r)) \, .
\end{equation}
For further  details including background material and the most recent developments centred around Lemma~DBC   see \cite{BeresnevichVelani7} and \cite[Section~2]{BHV2024}.

The upshot of Lemma~DBC is that in one way or another, quasi-independence on average is enough to guarantee full measure or equivalently that $
 A(x,N)  \to   \infty $ as $ N \to \infty
$ for $\mu$-almost all $x\in \Omega$. Clearly,  \eqref{fullind} implies \eqref{vbx1x} with $C=1$ and it is not surprising  that is  enough to imply the stronger asymptotic statement \eqref{asytriv}.  However, as with \eqref{fullind} when it comes to applications we rarely have \eqref{vbx1x} with $C=1$.  In the next section we discuss a weaker variant that  still enables us to  describe  the asymptotic behaviour of the counting function.

\subsection{A useful mechanism for establishing counting results} \label{mechCR}

The following  statement~\cite[Lemma~1.5]{Harman1998} represents an important  tool for establishing counting statements.  It has its bases in the familiar variance method of probability theory and   it can be viewed as the quantitative form of the divergence  Borel--Cantelli Lemma.  This will be made explicit in a moment.

\begin{lemma} \label{ebc}
Let $(\Omega,\mathcal{A},\mu)$ be a probability space, let $(f_n(x))_{n \in \N}$ be a sequence of non-negative $\mu$-measurable functions defined on $X$, and $(f_n)_{n \in \N },\ (\phi_n)_{n  \in \N}$ be sequences of real numbers  such that
$$ 0\leq f_n \leq \phi_n \hspace{7mm} (n=1,2,\ldots).  $$

\noindent 
Suppose that for arbitrary  $a,b \in \N$ with $a <  b$, we have
\begin{equation} \label{ebc_condition1}
\int_{\Omega} \left(\sum_{n=a}^{b} \big( f_n(x) -  f_n \big) \right)^2\mathrm{d}\mu(x)\, \leq\,  C\!\sum_{n=a}^{b}\phi_n
\end{equation}

\noindent for an absolute constant $C \ge 1$. Then, for any given $\varepsilon>0$,  we have
\begin{equation} \label{ebc_conclusion}
\sum_{n=1}^N f_n(x)\, =\, \sum_{n=1}^{N}f_n  \, +\, O\left(\Phi(N)^{1/2}\log^{\frac{3}{2}+\varepsilon}\Phi(N)+\max_{1\leq n\leq N}f_n\right)
\end{equation}
\noindent for $\mu$-almost all $x\in \Omega$, where $
\Phi(N):= \sum\limits_{n=1}^{N}\phi_n
$. 
\end{lemma}

Note that in statistical terms, if the sequence $f_n$ is the mean of $f_n(x)$; i.e.
$$
f_n = \int_{\Omega} f_n(x) \, \mathrm{d}\mu(x) \, ,
$$
then the l.h.s.~of~\eqref{ebc_condition1} is simply the variance ${\rm Var}(Z_{a,b}) $ of the random variable  $$Z_{a,b}=Z_{a,b}(x):=\sum_{n=a}^{b}f_n(x) \, . $$ In particular,
$
 {\rm Var}(Z_{a,b})  = \E(Z^2_{a,b})  -  \E(Z_{a,b})^2    \
{\rm where  }   \;
 \E(Z_{a,b}) = \int_{\Omega} Z_{a,b}(x) \, \mathrm{d}\mu(x)  \,
$  is the expectation of the random variable.

\medskip

We now use the above lemma to  explicitly give the quantitative form of Lemma~DBC involving the counting function \eqref{countdef}.   With this in mind, let  $A_n$ ($n \in \mathbb{N}$) be a family of  measurable subsets of $\Omega$   and consider Lemma \ref{ebc}  with
\begin{equation} \label{Harman_choice_parameters}
 \qquad f_n(x):= \one_{A_{n}}\!(x)  \qquad {\rm and } \qquad
 f_n:= \phi_n:= \mu(A_n) \, ,
\end{equation}
where   $  \one_{A_{n}}\!$   is the characteristic function  of the set  $A_n$ ($n \in \N$).     Then, clearly for any $ x \in \Omega$ and $N \in \N$  we have that
$$
{\rm l.h.s. \ of \ } \eqref{ebc_conclusion}  \ =   \ A(x, N)  \, :=  \,    \# \big\{ 1\le n \le   N :    x \in A_n  \big\}  \, ,
$$
where $A(x,N)$ is the counting function given by \eqref{countdef}.  Also,  observe that the main term on the ${\rm r.h.s. \ of }$ \eqref{ebc_conclusion} is
\begin{equation}  \label{needed}
\Phi(N) :=  \sum\limits_{n=1}^{N}\mu (A_n) \, . \end{equation}     Furthermore, it is easily verified  (see for instance \cite[Section~3]{levesley2024}) that for any $a,b \in \N$ with $a <  b$

%
%
%
%
%
%
%
\begin{eqnarray} \label{integral_square_estimate}
  {\rm l.h.s. \ of \ } \eqref{ebc_condition1}
   &= & \   \sum_{m,n=a}^b  \mu(A_m\cap A_n)  \  -  \  \left(\sum_{n=a}^{b}\mu(A_n) \right)^2\, .
\end{eqnarray}
\\  The following statement is now easily seen to be a direct consequence of Lemma~\ref{ebc}. \\

 \begin{lemmaQBC}
 \!
 \emph{ Let $(\Omega,\mathcal{A},\mu)$ be a probability space  and
let $(A_n)_{n \in \mathbb{N}}$ be a sequence   measurable subsets of $\Omega$.  Suppose there exists an absolute constant $C \ge 1$ such that
for arbitrary $a, b \in \N $ with $a < b$,
\begin{equation}\label{indie}
\sum_{n,m=a}^b  \mu(A_m\cap A_n)   \ \le \  \left(\sum_{n=a}^b  \mu(A_n)\right)^2  \ + \ C   \,  \sum_{n=a}^{b}  \mu(A_n)     \, .
\end{equation}
Then, for any given $\varepsilon>0$,  we have
\begin{equation} \label{ebc_conclusionqbc}
A(x, N)  \, =\, \Phi(N) \, +\, O\left(\Phi(N)^{1/2}\log^{\frac{3}{2}+\varepsilon}\Phi(N)  \right)
\end{equation}
\noindent for $\mu$-almost all $x\in \Omega$, where $
\Phi(N):= \sum\limits_{n=1}^{N}  \mu(A_n)
$.}
\end{lemmaQBC}

The upshot of Lemma~QBC is that if the sets $A_n$ are pairwise independent on average  (i.e., \eqref{vbx1x} with $C=1$) with an acceptable error term (i.e. as in \eqref{indie}) then we have an asymptotic statement with essentially the best error term  -- the square root of the main term is optimal (see for instance  the discussion following the statement of Lemma~1.5 in \cite{Harman1998}). 

\medskip

The results presented in this ``Borel--Cantelli'' section are for general $\limsup$ sets.  They should act as a guide to the type of measure and counting results we could potentially expect to prove for the specific  recurrent and shrinking target $\limsup$ sets - the main focus of this survey.

\section{Results for the shrinking target setup} \label{STSU}

We start with a quick recap of the problem under investigation.   With  \eqref{stques} in mind,  recall   there are two natural questions that fall under the general  umbrella of the ``shrinking target problem'', as  introduced in~\cite{hill1995}:
\begin{itemize}
  \item[]
  \begin{itemize} \item[\textbf{(P1)}] What is the $\mu$-measure of $W(T,\psi)$ if the measure sum in \eqref{stques} diverges? \\[-1ex]
  \item[\textbf{(P2)}] What is the Hausdorff dimension of $W(T,\psi)$  when as in \eqref{stques} the measure sum converges and so $\mu\big(W(T,\psi)\big)=0$?
      \end{itemize}
\end{itemize}
In \cite{hill1995} and the follow-up paper \cite{MR1471868},  the primary focus  was on the dynamics of expanding rational maps. Since then, the shrinking target problem has been studied in a broad array of dynamical settings. We refer the reader to \cite{aspenberg2019,barany2018,bugeaud2014,fang2020,ghosh2024,LI2023108994,li2014,persson2019} and the references therein for results concerning Hausdorff dimension and related fractal aspects, and to  \cite{allen2021,fernandez2012,galatolo2015,kim2007,tseng2008} and the references therein for measure-theoretic developments.

As mentioned in the introduction our primary  focus in this survey is on the measure aspect of the shrinking target problem. With this in mind, for both the shrinking target problem and the recurrent problem (where in the previous paragraph  \eqref{stques} and  $W(T,\psi)$ are  replaced by \eqref{rques}  and $R(T,\psi)$)   we will assume that the underlying   probability measure preserving system $(X,\mathcal{A},\mu,T)$ is not just ergodic but also exponentially mixing.  That is to say we assume that $T$ is     \emph{exponentially mixing\,} with respect to $\mu$ and so there exist constants $C> 0 $ and  $\gamma\in (0,1)$   such that
\begin{equation}\label{expmixballs}
  \Big|  \mu(B\cap T^{-n}F)   -   \mu(F)\mu(B) \Big|  \; \leq \;  C  \, \gamma^n   \mu(F)      \qquad  \forall  \ n \in \N \, ,
\end{equation}
for all balls $B$ in  $X$ and $\mu$-measurable sets $F$ in $X$.
The fact that we impose a condition beyond ergodic on the measure preserving system is necessary (see Example~\ref{needmore} below).

\begin{remark} \label{expind}
 Note that since $T$ is measure preserving,  by definition we have  that
 $$  \mu(T^{-n}F)     =  \mu(F)  \qquad \forall  \ F \in \mathcal{A} \, . $$  Thus, loosely speaking,   exponentially mixing  tells us that balls in $X$  and pre-images of measurable sets are pairwise independent (see \eqref{fullind}) up to an error term that is  exponentially decaying.
 \end{remark}

\begin{remark}
 Note that in certain situations, such as when $X$ is a subset of $\R^d$ and $\mu$ is a Radon measure,  the above notion of exponentially mixing  suffices  to imply \emph{strong-mixing}; that is
\[
\lim_{n\to \infty} \mu(E\cap T^{-n}F)=\mu(E)\mu(F)   \qquad \forall  \  E,F \in\cA \, .
\]
This in turn implies  ergodicity.
\end{remark}

\noindent  The reason for imposing  the condition that the system is exponentially mixing from the onset is two fold. The  main reason is that,  to the best of our knowledge, exponentially mixing in one form or another underpins all know results related to the recurrent setup  -- it is the basic assumption. The other  reason to impose the condition now, is that within the shrinking target setup we are able to clearly demonstrate its power. In short, as we shall soon see, exponential mixing gives us, essentially for free, the best possible ``expected'' outcome to the problem of determining what happens when the measure sum in \eqref{stques} diverges.  This is so unlike the recurrence setup (discussed in the next section) in which we need to  impose additional hypotheses and still have to work significantly if not infinitely harder!

For convenience,   we recall from the introduction  that the shrinking target set $W( T, \psi)$    is a $\limsup$ set. Indeed,
$
W( T, \psi) =  \limsup_{n \to \infty} W_n  \, ,
$
where for $n\in \N $\\
\begin{equation}  \label{souseful}
 W_n \,  =   \, W_n( T,\psi)
   \, := \,  \big\{ x \in X :  T^nx \in B_n:=B(x_0, \psi(n)) \big\}  \  =  \   T^{-n} (B_n) \, .
\end{equation}\\
Thus we are interested in points $x \in X$ whose orbit under $T$ `hits'  the target balls $B_n$ an infinite number of times.  With this in mind,  for $ x \in X$  and $N\in \N $, consider the counting function
\begin{eqnarray}
W(x,N;T,\psi)& :=  & \# \big\{ 1\le n \le   N :    x \in W_n \big\} \, \\[1ex]  & = & \# \big\{ 1\le n \le   N :    T^nx \in B_n  \big\}   \, .  \nonumber
   \end{eqnarray}
   We now show in a couple of lines, that if we have exponentially mixing then the sets $W_n$ satisfy \eqref{indie} and so in turn we are able to determine the asymptotic behaviour of the above  counting function  for $\mu$-almost all $x \in X$.

\smallskip Let  $m < n  $.  Then in view of  \eqref{souseful} and the fact that $T$ is measure preserving and exponentially mixing  with respect to $\mu$, it follows that there exist constants $C > 0 $ and $\gamma \in (0, 1) $ such that
\begin{eqnarray*} \label{rhssveasy}
\mu(W_m \cap W_n)  &  =  &  \mu\Big(T^{-m}(B_m)  \cap T^{-n}(B_n)\Big)  =  \mu\Big(B_m  \cap T^{-(n-m)}(B_n)\Big) \\[1ex]
& \le  &  \mu\big(B_m \big)  \mu\big(B_n \big)   +   C  \gamma^{n-m} \mu\big(B_n \big)  \ = \  \mu\big(W_m \big)  \mu\big(W_n \big)   +   C  \gamma^{n-m} \mu\big(W_n \big)  \, .
\end{eqnarray*}
In turn, it follows that for  arbitrary $a, b \in \N $ with $a < b$,
\begin{equation}\label{indiest}
\sum_{n,m=a}^b  \mu(W_m\cap W_n)   \ \le \  \left(\sum_{n=a}^b  \mu(W_n)\right)^2  \ + \   \sum_{n=a}^{b}  \mu(W_n)  \, \Big( 1 + C     \, \sum_{m=1}^{\infty}  \gamma^m  \Big)    .
\end{equation}
Now the sum involving $\gamma$ is convergent and so this shows that the sets $W_n$ satisfy the independence condition  \eqref{indie}. Thus Lemma~QBC implies the following  statement.   It  tells us that if the `natural'  measure sum  diverges  then for  $\mu$-almost all $x \in X$  the orbit  `hits' the target sets $B_n$ the `expected' number of times.

\begin{theorem}\label{shrinktarg}
Let $(X,\mathcal{A},\mu,T)$ be a measure-preserving dynamical system  and  suppose that $T$ is exponentially mixing with respect to $\mu$.  Let $\psi:\N\to\R_{\ge 0}$ be a  real, positive function.  Then, for any given $\varepsilon>0$, we have that
\begin{equation} \label{shrinktargcount}
W(x,N; T, \psi):= \sum_{n=1}^N\one_{B_n}(T^nx) =   \Phi(N)+O\left(\Phi^{1/2}(N) \ (\log\Phi(N))^{3/2+\varepsilon}\right)
\end{equation}
\noindent for $\mu$-almost all $x\in X$, where
\begin{equation}
\label{shrinktargsum}
\Phi(N):=\sum_{n=1}^N \mu(B_n) \,  .
\end{equation}
\end{theorem}

\medskip

 \begin{remark} Within the context of the above theorem, it is in fact  possible to replace the exponentially mixing assumption by the weaker notion of \s-mixing (short for summable mixing) -- see \cite[\S1.1]{LI2023108994}  for the more general statement and for connections to other works.  On a different note, and keeping in mind the comment made immediately after the statement of Lemma~QBC, we observe that the error term in Theorem~\ref{shrinktarg} is essentially best possible.\end{remark}

\bigskip

Note that  the measure sum \eqref{shrinktargsum} is equivalent to
\begin{equation} \label{shrinktargsum2}  \Psi(N):=\sum_{n=1}^N\mu\big(W_n\big)   \, . \end{equation}
Thus,  the theorem shows that for $\mu$-almost all $x \in X $,  the asymptotic behaviour  of the counting function $W(x,N;T,\psi)$  is determined by the behaviour of   the measure sum  $ \Psi(N)$ involving the ``building block''  sets $W_n$ associated with the $\limsup$ set  $W(T,\psi)$.   This together with the fact that  $\Psi(N)$ is independent of $x \in X $,
is well worth keeping in mind for future  comparison with the analogous recurrent  problem.     Next note that  by definition, $ x \in W(T,\psi) $ if and only if $\lim_{N \to \infty}  W(x,N;T,\psi) = \infty$ and  so  an immediate consequence of Theorem~\ref{shrinktarg} is the following zero-one  measure criterion (which naturally is in line with the convergent and divergent Borel--Cantelli Lemmas for  general $\limsup$ sets).

\begin{theorem} \label{shrinktargcor}
	Let $(X,\mathcal{A},\mu,T)$ be a measure-preserving dynamical system  and  suppose that $T$ is exponentially mixing with respect to $\mu$.   Let $\psi:\N\to\R_{\ge 0}$ be a  real, positive function.  Then
\begin{eqnarray}\label{appdiv}
		\mu\big(W(T,\psi)\big)=
		\begin{cases}
			0 &\text{if}\ \  \sum_{n=1}^\infty \mu\big(B_n\big)<\infty\\[2ex]
			1 &\text{if}\ \ \sum_{n=1}^\infty \mu\big(B_n\big)=\infty.
		\end{cases}
	\end{eqnarray}
\end{theorem}

\medskip

We conclude this section with an example of a classical dynamical system that illustrates a key point: to establish such a zero-one measure criterion (let alone a quantitative statement) in the shrinking target setup, it is necessary to assume more than mere ergodicity and measure preservation.

\begin{example}  \label{needmore}
 Let $X=\R/\Z$  and let $T$ be the irrational rotation by $\alpha\in\R\setminus\Q$ defined as
\[T:X\to X: x \to T(x):=\alpha+x\quad(\rm{mod}\ 1).\]
Furthermore, let $\mu$ be the one-dimensional Lebesgue measure restricted to $[0,1]$.  Then it is well known that the corresponding  dynamical system  $(X,\mathcal{A},\mu,T)$ is  measure preserving  and ergodic but not mixing.
In order to show that this system gives rise to the desired  shrinking target counterexample, a little background from the theory of Diophantine approximation is required.  Throughout, for $ x \in \R$, let
$
  \|x\| := \min \{ | x-m | : m \in \Z \}
$
denote the distance from $x$ to the nearest integer.   A classical result of Dirichlet (1842)  states that for any $x \in \R$
$$
\| n x \|   \, \le \, n^{-1}    \  \ \  {\rm for   \ infinitely \  many \  \ }   n \in \N \, .  
$$
Now, with this in mind, the set ${\rm Bad}$ consists of those real numbers for which Dirichlet's theorem cannot  be improved by an arbitrary constant.  Not surprisingly, such numbers are  referred to as \emph{badly approximable numbers} and formally $ x \in {\rm Bad}$ if
$$
\liminf_{n \to \infty}  n \, \| n x\| >0    \, .
$$
 For the sake of completeness, we mention that ${\rm Bad}$ is a set of Lebesgue measure zero but full Hausdorff dimension; that is to say $\dim  \, {\rm Bad}  =1$.    For background and  further information  see  \cite{RothVol,BRV2016} and references within.

 \medskip

 \noindent \emph{Shrinking target counterexample. }   Let $\psi:\N\to\R_{\ge 0}$ be a  real, positive decreasing function such that $\psi(n) \to 0 $ as $n \to \infty$.  By definition:
 \begin{eqnarray}   \label{egshr}
 x \in  W(T, \psi) & \quad  \Longleftrightarrow  \quad &  T^nx \in B(x_0, \psi(n))  \ \  \  {\rm for  \ \ i.  \ m. \ }   n \in \N   \nonumber \\[1ex]
 & \quad  \Longleftrightarrow  \quad & \|  n \alpha + x - x_0    \|    \le  \psi(n)   \  \ \  {\rm for  \ \ i.  \ m. \ }   n \in \N
 \end{eqnarray}
 Thus, the shrinking target set corresponds to the so called `` twisted'' inhomogeneous approximation set
 $
 \mathcal{T}_{\alpha}(\psi)
 $  consisting of real numbers $s=x-x_0\in X$ such that
$$
\| n \alpha - s\|   \, \le \, \psi(n)
$$
for infinitely many $n \in \N$. A beautiful result of Kurzweil \cite{Kurzweil2} states that for a
$$
\mu(  \mathcal{T}_{\alpha}(\psi) )  = 1   \quad \forall  \  \  \psi  \in  \mathcal{D}  \quad  \Longleftrightarrow  \quad \alpha \in {\rm Bad} \, ,
$$
where $  \mathcal{D} $ is the set all real, positive decreasing functions $\psi$ such that $\sum_{n=1}^{\infty} \psi(n) = \infty $.  For background and  further information  regarding the ``twisted'' theory of Diophantine approximation see  \cite[Section~9]{RothVol} and references within.  Anyway, the upshot of Kurzweil's Theorem is that for any $ \alpha \notin  {\rm Bad}$, there exists a $\psi$ such that $\sum_{n=1}^{\infty} \psi(n) = \infty $ but  $\mu(\mathcal{T}_{\alpha}(\psi)) \neq 1 $.  Now for any $n \in \N$ , it is easily verified that $\mu(W_n) \asymp  \psi(n) $ and so it follows that 
$$\sum_{n=1}^{\infty} \mu( W_n)    \asymp    \sum_{n=1}^{\infty} \psi(n) = \infty   \quad  { \rm but }  \quad   \mu \big(W(T, \psi))   \neq 1   \, . $$
Here of course  the sets $W_n$ are the  ``building block''  sets  (see  \eqref{souseful}) associated with the $\limsup$ set  $W(T,\psi)$.

\medskip

It is interesting to note that the classical ``circle rotation'' dynamical system does not lead to a recurrence counterexample.  Let $\psi:\N\to\R_{\ge 0}$ be a  real, positive function.  By definition:
 \begin{eqnarray}   \label{egrec}
 x \in  R(T, \psi) & \quad  \Longleftrightarrow  \quad &  T^nx \in B(x, \psi(n))  \ \ \  {\rm for \  \ i.  \ m. \ }   n \in \N   \nonumber \\[1ex]
 & \quad  \Longleftrightarrow  \quad & \|  n \alpha + x - x    \|    \le  \psi(n)    \ \ \  {\rm for  \ \ i.  \ m. \ }   n \in \N    \nonumber \\[1ex]
 & \quad  \Longleftrightarrow  \quad & \|  n \alpha    \|    \le  \psi(n)    \ \ \  {\rm for \ \ i.  \ m. \ }   n \in \N  \nonumber
 \end{eqnarray}
 In particular, it follows that 
 \begin{eqnarray*}
		R_n=
		\begin{cases}
			X &\text{if}\ \  \|  n \alpha    \|    \le  \psi(n)   \\[2ex]
			\emptyset &\text{if}\  \ \|  n \alpha    \|    >  \psi(n)  \, 
		\end{cases}
	\end{eqnarray*}
 and so
 \begin{eqnarray}  \label{poil}
		R(T,\psi)=
		\begin{cases}
			X &\text{if}\ \  \sum_{n=1}^\infty \mu\big(R_n\big)=\infty\\[2ex]
			\emptyset &\text{if}\ \ \sum_{n=1}^\infty \mu\big(R_n\big)<\infty.
		\end{cases}
	\end{eqnarray}  
    
\noindent Here of course  the sets $R_n$ are the  ``building block''  sets  (see  \eqref{def_An}) associated with the $\limsup$ set  $R(T,\psi)$.
 Clearly \eqref{poil}, implies that $\mu\big(R(T,\psi)\big)= 1$ (resp. 0) if $\sum_{n=1}^\infty \mu\big(R_n\big)$ diverges  (resp. converges) which is perfectly in line the ``expected'' zero-one measure criterion. 
\end{example}

\begin{remark}
To the best of our knowledge, there is no known simple example in the recurrence setting demonstrating that, if $(X,\mathcal{A},\mu,T)$  is an ergodic, measure-preserving dynamical system and $\mu$ is a  uniform measure (e.g., Ahlfors regular), then the ``expected'' zero-one measure criterion fails\footnote{We would like to thank Ayesha Bennett and Niko Nikov for  bringing this  ``hole'' in the literature to our  attention.} (cf. Example ABB in the next section, where $\mu$ is non-uniform).  Recall, a measure $\mu$ on a metric space $(X,d)$ is \emph{$\tau$-Ahlfors regular} if there exists a constant  $C\ge 1$  such that for any ball $B(x, r)  \subset X$ with $x \in X$
\begin{equation}  \label{tar}
C^{-1}r^{\tau}\leq \mu(B(x,r))\leq Cr^{\tau}
\,\qquad  \forall \ 0<r\leq |X|
\, ,  \medskip
\end{equation}
where $|X|$ denotes the diameter of $X$.
\end{remark}

%
%
%
%

\section{Results for the recurrence setup} \label{RSU}

For convenience,   we recall from the introduction  that the recurrence  set $R( T, \psi)$    is a $\limsup$ set. Indeed,
$
R( T, \psi) =  \limsup_{n \to \infty} R_n  \, ,
$
where for $n\in \N $
\begin{equation}  \label{anrpsidef}  
 R_n \,  =   \, R_n( T,\psi)  \,  :=  \,  \big\{ x \in X :  T^nx \in B(x, \psi(n)) \big\}
  \, .
\end{equation}
As mentioned in the introduction, under the very basic assumption that $(X,\mathcal{A},\mu,T)$  is a measure-preserving dynamical system, we know that $ 
  \mu( R(T,c) )   = 1 $ if $\psi$ is a constant $c$ and so we are interested in determining the ``size'' of  $R(T,\psi)$ when  $\psi(n) \to 0 $ as $n \to \infty$. With this in mind, the first results date back to the pioneering work of  Boshernitzan \cite{boshernitzan1993} who studied the case $\psi(n)=n^{-1/\alpha}$ ($\alpha>0$).  More precisely, he showed that if the $\alpha$-dimensional Hausdorff measure $\cH^{\alpha}$ of   $X$ is zero for some $\alpha > 0$, then for any $\epsilon > 0$
  $$
  \mu( R(T,\epsilon \psi_{\alpha} ))   = 1   \qquad {\rm where } \qquad \psi_{\alpha}: n \to n^{-1/\alpha}    \, .    
  $$ 
  On the other hand, if  $\cH^{\alpha}(X) > 0$, then for $\mu$-almost all $ x \in X$, there is a constant $c(x) > 0 $  such that 
  $$
   T^nx \in B\big(x,c(x)  \psi_{\alpha}(n) \big) \  \ \  {\rm for  \ \ i.  \ m. \ }   n \in \N  \, .
  $$
  Subsequently, 
  under various additional assumptions on the  measure-preserving dynamical system, analogues  of the  shrinking target zero-one measure criterion  given by Theorem~\ref{shrinktargcor} have been established for the recurrent set $R(T,\psi) $ in numerous works  -- see for instance  \cite{baker2021,baker2024,chang2019,he2022,hussain2022,Kirsebom2023,kleinbock2023} and the references within.    We refer the reader to \cite{he2022,hu2024,seuret2015,tan2011,wu2025} and the references within for results concerning Hausdorff dimension and related fractal aspects.   
  
  Focusing on the measure-theoretic rather than dimension developments, we note that Boshernitzan's result requires only that the system in question is measure preserving. 
To the best of our knowledge,  this condition alone is insufficient to derive  a zero-one measure criterion.  Indeed, although  some form of exponential mixing underpins all  of the aforementioned  results in the recurrent setting, and this appears to be the fundamental assumption, additional hypotheses are typically  necessary.  In particular, beyond exponential mixing, common supplementary conditions imposed on the measure-preserving system include properties such as expansiveness, bounded distortion, conformality, and Ahlfors regularity.  Some combination of these “extra” properties - whether implicitly or explicitly -appears in all the works cited above. For precise definitions, we refer the reader to \cite[Section~1]{hussain2022}, where each of these conditions is required and stated explicitly.

We stress that none of the ``extra'' properties mentioned above are required  within the shrinking target setup -- even when establishing  the stronger quantitative statement of Theorem~\ref{shrinktarg}. So it is natural to wonder why they are necessary in the recurrence setup.   In short, it is because the  ``building block''  sets  $R_n$ (see  \eqref{anrpsidef}) associated with the $\limsup$ set  $R(T,\psi)$ are not  the pre-image of  balls.   This means that we cannot directly  exploit ``measure preserving'' and ``exponential mixing''  to show that the sets $R_n$ ($n\in \N$) are quasi-independence on average (see \eqref{vbx1x}),  let alone satisfy the stronger independence condition \eqref{indie}, as is the case in the shrinking target setup (see \eqref{indiest}).     It turns out that, in order to exploit exponential mixing, one must work locally, since locally the set   $ R_n$ can be expressed as the pre-image of a ball. More precisely, at the most basic level, the strategy hinges on the following observation, which is a straightforward application of the triangle inequality. For the details see for instance \cite[Lemma~2.2]{hussain2022}.

\begin{lemma}
Let $B=B(z, r)$ be a ball centred  at $z \in X$ and radius $r > 0$.  Then for any $m \in \N$  with  $\psi(m) > r$
    \begin{equation*}
      B \cap  T^{-m} \big( B \left(z,  \psi(m) - r \right) \big)   \ \subset \   B \cap R_m     \ \subset \   B \cap  T^{-m} \big( B \left(z,  \psi(m) + r  \right) \big)  \, . 
      \end{equation*}
\end{lemma}


\noindent Effectively implementing this local decomposition of $R_n$ to subsequently leverage exponential mixing is - loosely speaking - precisely where the additional assumptions (such as expanding, bounded distortion, conformality,
and Ahlfors regularity)  come into play.

The works mentioned above  (namely, \cite{baker2021,baker2024,chang2019,he2022,hussain2022,Kirsebom2023,kleinbock2023})  are all centred around establishing a zero-one measure criterion for the particular class of measure preserving system under consideration.  In short, the associated recurrent sets  are shown to satisfy the ``expected'' zero-one measure criterion (see \eqref{divcondsv1} below).  The  ``conformal'' measure preserving  dynamical systems considered in \cite{hussain2022} is among the most general studied  to date, and includes many well-known systems that admit a “nice” countable partition - such as the Gauss map.  In contrast,  
the results in the earlier works \cite{baker2021,chang2019} apply to finite conformal iterated function systems (with the open set condition) and thus do not cover the Gauss map or even the L\"{u}roth map -- a simpler linear analogue  of the Gauss map. This completes our brief overview of known systems for which the  ``expected'' zero-one measure criterion holds.

 We now turn our attention to the problem of determining the stronger quantitative  statement, which, within the shrinking target setup (namely   Theorem~\ref{shrinktarg}),  we essentially obtained for free.   Things now get interesting.  Even for finite conformal iterated functions systems for which the  ``expected'' zero-one measure criterion holds, the corresponding ``expected'' asymptotic statement is not generally true!   Furthermore,  we shall see that if one considers  non-Ahlfors regular measures associated with finite conformal iterated functions systems, then even the ``expected'' zero-one measure criterion can fail!

At this point, we will closely follow the exposition presented  in \cite[Section 1.3]{jungie}  where the measure preserving system  is  a self-conformal system on $\R^d$.  We  also stress that all the results appearing below are  established \cite{jungie}.  As in \cite{jungie}, for a self-conformal system on~$\mathbb{R}^d$, we adopt the notation~$(\Phi, K, \mu, T)$ in place of~$(X, \mathcal{A}, \mu, T)$. Briefly, the setup is as follows:
  \begin{itemize}
    \item  $\Phi=\{\varphi_j\}_{1\leq j\leq m}$ ($m\geq2$) is a $C^{1+\alpha}$ conformal IFS (iterated function system)  on $\R^d$ satisfying the open set condition.
    \vspace*{1ex}
    \item $K\subseteq\R^d$ is the self-conformal set generated by $\Phi$.
      \vspace*{1ex}
    \item $\mu$ is a  Gibbs measure on $K$.
      \vspace*{1ex}
    \item $T:\R^d\to\R^d$  is a natural map induced by $\Phi$ such that $T|_K:K\to K$ is conjugate to the shift map on the symbolic space $\{1,2,...,m\}^{\N}$.
\end{itemize}

The precise definitions of each of the above properties can be found in \cite[Section~2]{jungie}. For the purposes of this survey, however, the examples given below should suffice. 
Nevertheless, it is worth noting that, by definition, any self-conformal system naturally inherits the properties of expanding, bounded distortion, and conformality - each of which features in the works cited above dealing with conformal dynamical systems. However, the associated Gibbs measure  $\mu$  need not be Ahlfors regular, in contrast to those works.  Until recently, for $d \ge 2$ it was not known whether self-conformal systems exhibit exponential mixing.  The main result in \cite[Theorem~1.2]{jungie}  implies that if $(\Phi,K,\mu,T)$ is  a self-conformal system on $\R^d$ then $T$ is exponentially mixing with respect to $\mu$.  
For $d \ge 2$, this  appears to be the first result demonstrating  the existence of a natural class of  measure preserving systems whose associated  maps are exponentially mixing with respect to non-trivial fractal measures.

In any case the upshot of main result in \cite{jungie}  is that  Theorem~\ref{shrinktarg} holds for any  self-conformal system on $\R^d$.  We now turn to the task of determining its analogue within the recurrence setup. With \eqref{anrpsidef} in mind, given $ N \in \N$ and $\x \in K$,  consider  the  counting function
\begin{eqnarray} \label{countdefst}
R(\x,N;T,\psi)& :=  & \# \big\{ 1\le n \le   N :    \x \in R_n(T,\psi) \big\} \, \nonumber  \\[1ex]  & = & \# \big\{ 1\le n \le   N :    T^n\x \in B(\x , \psi(n))  \big\} \nonumber \\[1ex] & = &  \textstyle{\sum_{n=1}^{N}}\mathbbm{1}_{B(\x,\psi(n))}(T^n\x)  \, .
   \end{eqnarray}
   In line with the shrinking target setup - and more broadly, the quantitative Borel–Cantelli framework (namely, Lemma~QBC) - it is natural to expect that, for  $\mu$-almost all $\x \in K$,    the asymptotic behaviour  of this counting function  is governed  by the $\mu$-measure sum  of the sets $R_n:=R_n(T,\psi)$.  More precisely, the following statement represents the ``expected'' asymptotic statement.

\begin{claimf}  \label{claimf}
  \emph{ Let $(\Phi,K,\mu,T)$ be a self-conformal system on $\R^d$, let $\psi:\N\to\R_{\geq0}$ be a real positive function and assume that $\sum_{n=1}^{\infty}\mu(R_n(T,\psi)) $ diverges.  Then, for $\mu$-almost  all $\x \in  K$ 
 \begin{equation}\label{sbc}
 \lim_{N\to\infty}\frac{\sum_{n=1}^{N}
 \mathbbm{1}_{B(\x,\psi(n))}(T^n\x)}{\sum_{n=1}^{N}\mu(R_n(T,\psi))} \, = \, 1 \, .
 \end{equation}}
 \end{claimf}

 \noindent Such a claim was also  alluded to in \cite[Section~1]{levesley2024} and it was  shown to be true for a large class of piecewise linear maps in $\R^d$. However, it turns out that in general the claim is  false (hence the label ``F'') in a rather strong sense.  Indeed,  we are able to give   explicit examples  of self-conformal systems for which the $\mu$-measure of the $\limsup$ set $R(T,\psi)$ is one but the limit
appearing in \eqref{sbc}  is not even a constant let alone one (cf. Example~ABB below). In
other words, even after excluding a set of $\mu$-measure zero, the limit in \eqref{sbc}  depends on $\x$
and thus for these self-conformal systems the associated recurrent sets exhibit (unexpected
and extreme) behaviour that is not present for shrinking target sets.   The following summarises the counterexamples to the claim.  The full details  are  given in \cite[Section~7.4]{jungie}.

   \begin{example}  \label{egcantor}
Let
 $\Phi=\{\varphi_1,\varphi_2\}$ where
   $\varphi_1:[0,1]\to[0,1/3]$ and $\varphi_2:[0,1]\to[2/3,1]$ are  given by
    \[
    \varphi_1(x)=\frac{x}{3},  \qquad \varphi_2(x)=\frac{x+2}{3}  \qquad \forall \  x\in[0,1]  .
    \]
 Then  $\Phi$ is the well-known  conformal IFS with the attractor $K$ being the standard middle-third Cantor set.  It gives rise to the ``natural'' self conformal system  $(\Phi,K,\mu,T)$ in which  $\mu:=\cH^{\tau}|_K$  ($\tau=\log2/\log 3$) is the standard Cantor measure and  $T:[0,1]\to[0,1]$ is given by $ Tx=3x~\tmod~1   \, . $
   Now let $\psi:\N\to\R_{\geq0}$ be the constant function  given by $\psi(n) \, := \,  \textstyle{ \Large{\frac{1}{3}+\frac{2}{3^2} }}  \, . $
    Then, for $\mu$-almost all  $x\in K$, we have that 
    
   \begin{eqnarray*} \label{klj}
      \lim_{N\to\infty}\frac{\sum_{n=1}^{N}
 \mathbbm{1}_{B(x,\psi(n))}(T^nx)}{\sum_{n=1}^{N}\mu(R_n(T,\psi))}
      =\left\{
\begin{aligned}
    &\textstyle{\frac{4}{5}  \quad \text{if}  \quad x\in\left(\big[0,\frac{1}{9}\big]\cup\big[\frac{8}{9},1\big]\right)\cap K} ,\\[2ex]
    &\textstyle{\frac{6}{5}   \quad \text{if} \quad x\in\left(\big[\frac{2}{9},\frac{1}{3}\big]\cup\big[\frac{2}{3},\frac{7}{9}\big]\right)\cap K}.
\end{aligned}
      \right.
   \end{eqnarray*}
\end{example}

\bigskip

 \begin{example}  \label{egfull}
Let $\Phi=\{\varphi_1,\varphi_2,\varphi_3,\varphi_4\}$ be a conformal IFS on $[0,1]$ given by
 \begin{eqnarray*}
     \varphi_1(x)=\frac{1}{4}x,\ \  \varphi_2(x)=\frac{1}{2(1+x)},\ \ \varphi_3(x)=\frac{1+x}{2+x},\ \  \varphi_4(x)=\frac{2}{2+x} \qquad \forall \  x\in[0,1]    \,.
 \end{eqnarray*}
    Then  it can be  verified that $\Phi$ gives rise to a  self-conformal system $(\Phi,K,\mu,T)$ in which $K$ is  the unit interval,  $\mu$ is the natural Gibbs measure supports on $K$ that is absolutely continuous with respect to Lebsegue measure and $T:[0,1]\to[0,1]$ is given by
    \begin{eqnarray*}
        Tx=\left\{
        \begin{aligned}
            &  \textstyle{ 4x,  \qquad \  \text{if}\ \ 0\leq x<\frac{1}{4},}\\
            & \textstyle{ \frac{1}{2x}-1,~~~~ \ \text{if}\ \ \frac{1}{4}\leq x< \frac{1}{2},}\\[4pt]
            & \textstyle{\frac{2x-1}{1-x}, \ \ \ \ ~~~~~~\text{if}\ \ \frac{1}{2}\leq x<\frac{2}{3},}\\[4pt]
            &\textstyle{\frac{2}{x}-2, \quad \text{if}\ \ \frac{2}{3}\leq x\leq 1.}
        \end{aligned}\right.
    \end{eqnarray*}
Now let    $\psi:\N\to\R_{\geq0}$  be a real positive function such that $\psi(n) \to 0 $ as $n \to \infty$ and  $\sum_{n=1}^{\infty}\mu(R_n(T,\psi)) = \infty$. Then, for $\mu$-almost all  $x\in K:=[0,1]$, we have that
\begin{equation*}
     \lim_{N\to\infty}\frac{\sum_{n=1}^{N}
 \mathbbm{1}_{B(x,\psi(n))}(T^nx)}{\sum_{n=1}^{N}\mu(R_n(T,\psi))} \, = \,
     \frac{2\log 2}{1+x}.
\end{equation*}
\end{example}

\bigskip

 While the first more familiar ``Cantor''  example requires less sophisticated tools to setup and execute,  it does rely on $ \psi$ being a constant function.  It is also worth pointing out that that even though Example~\ref{egcantor} corresponds to a piecewise linear system it is not covered by \cite{levesley2024} since the  measure $\mu$ is not Lebesgue measure.

\medskip

\begin{remark}
    The counterexamples show that even though we have exponentially mixing for self-conformal systems (see \cite[Theorem 1.1]{jungie} for the details) we can not in general guarantee  that the sets $R_n(T,\psi)$ satisfy the strong quasi-independence average  condition \eqref{indie}
    as in the shrinking target framework. The point is that if they did then Lemma~QBC (quantitative  Borel--Cantelli) would  establish Claim~F.
\end{remark}

\medskip

Note that in both Example~\ref{egcantor} and \ref{egfull},  we still have that
 $\lim_{N \to \infty} R(x,N;T,\psi) = \infty  $  for $\mu$-almost all $x \in K$  and so  $\mu(R(T,\psi))=1$.  Moreover, we  highlight  the fact that in both examples  the measure $\mu$ is Ahlfors regular  and that for such measures this phenomena (under the assumption that $\sum_{n=1}^{\infty}\mu(R_n(T,\psi)) $ diverges) is known to hold for any self-conformal system  (see \cite{baker2021})
  and indeed for the  more general systems considered in  \cite{hussain2022}.  
The upshot of the above is that given a self-conformal system $(\Phi,K,\mu,T)$ on $\R^d$ for which the Gibbs measure $\mu$ is Ahlfors regular, and a real positive function  $\psi:\N\to\R_{\geq0}$,  then we obtain the ``expected'' zero-one measure criterion:
\begin{eqnarray}  \label{divcondsv1}
		 \mu\left(R(T,\psi)\right)=
		\begin{cases}
			0 &\text{if}\ \  \sum_{n=1}^{\infty}\mu\big(R_n(T,\psi)\big)<\infty,\\[2ex]
			1 &\text{if}\ \ \sum_{n=1}^{\infty}\mu\big(R_n(T,\psi) \big)=\infty.
		\end{cases}
\end{eqnarray}\\

\begin{remark}  \label{HMcoin} It is easily verified that within the setup of self-conformal systems, the notion of  $\mu$ being equivalent to $\cH^{\tau}|_K$  and $\mu$ being $\tau$-Ahlfors regular (see \eqref{tar}) coincide  -- for the details, see the proof of Theorem~2.7 in \cite{fan1999}. In general, we only have that the latter implies  the former. Recall,  that we say  Borel measures $\mu$ and $\nu$ on a metric space $(X,d)$ are equivalent if $\nu(E)\asymp \mu(E) $   for any Borel subset $ E\subseteq X$.
  \end{remark}

 \noindent 
 As we have already seen the convergent part of  \eqref{divcondsv1} is  a straightforward consequence of the  convergent Borel--Cantelli Lemma (see \S\ref{theBCsec}, Lemma~CBC).    In view of this, it is tempting to
    suspect that  at the coarser level of a zero-one measure criterion the analogue of  Claim~F is true; that is to say that \eqref{divcondsv1} is true for any self-conformal system on $\R^d$.  
     However, this  turns out not to be the case.  In a recent  beautiful paper,  Allen, Baker $\&$ B\'{a}r\'{a}ny \cite{allen2025}  consider the recurrent problem within the symbolic dynamics setting for topologically mixing sub-shifts of finite type. More precisely, in this setting they provide sufficient conditions for $\mu(R(T,\psi))$ to be zero or one when $\mu$ is assumed to be a non-uniform Gibbs measure and thus is not Ahlfors regular.  In terms of Bernoulli measures defined on the full shift,  the condition on the measure means that the  components of the defining probability vector  are not all equal.    As a consequence of their main result, 
     when specialised to the middle-third  Cantor we obtain the following concrete example that shows that \eqref{divcondsv1} is not true for any self-conformal system.

\medskip

  \noindent \textbf{Example~ABB.} \ Let $\Phi=\{\varphi_1 \, , \varphi_2 \}$, $T$ and  $K$ be as in Example~\ref{egcantor}.  Recall, $K$ is the  standard  middle-third Cantor. Now let $\mu$ be the weighted Cantor measure  associated with the probability vector  $(p_1,p_2)$ with $p_1\neq p_2$.  Let $\alpha>0$ and $\psi_{\alpha}(n)=3^{-\linte{\alpha\log n}}$. If
    \[
    \frac{1}{-(p_1\log p_1+p_2\log p_2)}<\alpha<\frac{1}{-\log(p_1^2+p_2^2)},
    \]
    then
    \[
    \sum_{n=1}^{\infty}\mu(R_n(T,\psi_{\alpha}))=\infty\qquad\text{but}\qquad \mu(R(T,\psi_{\alpha}))=0.
    \]

\medskip

A straightforward consequence of  Example~ABB is that for $\mu$-almost all $x \in K$
    \[
   \sum_{n=1}^{\infty}\one_{B(x,\psi_{\alpha}(n))}
(T^nx)  \, \ll \, 1   \quad {\rm and \ so \ }  \quad
    \lim_{N\to\infty}\frac{\sum_{n=1}^N\one_{B(x,\psi_{\alpha}(n))}
(T^nx)}{\sum_{n=1}^N\mu(R_n(T,\psi_{\alpha}))}=0  \, .
    \]
 In other words, even though the limit is a constant  for $\mu$-almost all $x \in K$,  it is not one (cf.  Claim~F). Note that  Examples \ref{egcantor} and \ref{egfull} show that Claim~F is false  even when  $\mu(R(T,\psi))=1$ and that for $\mu$-almost all $x \in K$, the limit under consideration is dependent on $x$ and thus  not a constant; that is to say that Claim~F is false on a large scale!


\medskip

Given that Claim~F is false, it is natural to  attempt to establish   an appropriate ``modified'' statement  that is true for the  full range of dynamical systems under consideration (namely, self-conformal systems).  Such a statement  would obviously follow on establishing the analogue  of Theorem~\ref{shrinktarg} for recurrent sets.  Indeed,  this is the ultimate goal as it  would provide an asymptotic result with an error term.    With this in mind, in order to state the first main result (for recurrent sets) we need to introduce a particular function that will determine the appropriate setup and thus the  asymptotic behaviour.   As usual,  let
$(\Phi,K,\mu,T)$   be a self-conformal system and a $\psi:\N\to\R_{\geq0}$ be a real,  positive function. Then for each  $n\in \N$,  we define the function
%
 $$
 t_n(\cdot)=t_n(\cdot,\psi):K\to\R_{\geq0}
 $$
  by
  \begin{eqnarray}\label{defofrn}
    t_n(\x)=t_n(\x,\psi):=\inf\left\{r\geq0:\mu(B(\x,r))\geq\psi(n)\right\}
\end{eqnarray}
if $ \psi(n)  \leq  1 $  and  we put $ t_n(\x)$ equal to the diameter of the bounded set  $K$ otherwise.
  With the definition of $t_n$ in mind,  and by  exploiting the fact that self-conformal  systems are exponentially mixing,  we are able to establish the following analogue  of Theorem~\ref{shrinktarg} for recurrent sets.   For the proof see \cite[Section~7.1]{jungie}.

\begin{theorem}\label{quantrec}
    Let $(\Phi,K,\mu,T)$ be a self-conformal system on $\R^d$ and let $\psi:\N\to\R_{\geq0}$ be a real positive function such that $\psi(n) \to 0 $ as $n \to \infty$. Furthermore, for $n \in \N$
    let $t_n :K\to\R_{\geq0}$ be given by  (\ref{defofrn}). Then for any $\epsilon>0$, we have
    \begin{equation} \label{recurrentcount}
        \sum_{n=1}^{N}\one_{B(\x,t_n(\x))}(T^n\x)=\Psi(N)+O\left(\Psi(N)^{1/2}\log^{\frac{3}{2}+\epsilon}(\Psi(N))\right)
    \end{equation}
    for $\mu$-almost all $\x\in K$, where
     \begin{equation} \label{recurrentsum}
    \displaystyle\Psi(N):=\sum_{n=1}^N\psi(n)   \, .
    \end{equation}
\end{theorem}

\noindent
\begin{remark} \label{remarkthm1.5}  Several comments are in order.
\begin{enumerate}
  \item[(i)] It turns out (see \cite[Lemma~7.3]{jungie}) that for all  $ \x \in K $ and all sufficiently large  $n \in \N$
\begin{equation}  \label{omg}   \mu(B(\x,t_n(\x))) = \psi(n)   \, .   \end{equation}
  Thus, up to an additive constant, the sum \eqref{recurrentsum} is simply the sum of the $\mu$-measure of the ``target balls'' $B(\x,t_n(\x)) $  associated with the modified  counting function appearing  on the left hand side of~\eqref{recurrentcount}.  In short, if the measure $\mu$ is non-uniform then the measure of a ball $B(\x, r)$  depends on its location $\x$  and not just its radius $r$.  In order to take this into account,  for  $n$ large, the radii of the target balls within the framework of  Theorem~\ref{quantrec}   are adjusted so that they all have the same measure  (namely $\psi(n)$) regardless of location.
  \medskip
  \item[(ii)] Let $ \hat{R}_n(\x,N;T,\psi) $   denote the  modified  counting function appearing  on the left hand side of~\eqref{recurrentcount}.  Then by definition,
  \begin{equation*} \label{countdefst-hat}
\hat{R}_n(\x,N;T,\psi)  =  \# \big\{ 1\le n \le   N :    \x \in  \hat{R}_n(T,\psi) \big\} \, ,
   \end{equation*}
   where
   $$
  \hat{R}_n(T,\psi)   :=    \big\{\x\in K :    T^n\x \in B(\x , t_n(\x))  \big\} \, .
   $$
 It turns out  (see \cite[Lemma~7.5]{jungie}) that there exists a constant $0< \gamma < 1 $ such that
$$ \mu( \hat{R}_n(T,\psi))=\psi(n)+O(\gamma^n)    \, . $$
The upshot of this and the equality \eqref{omg} appearing in (i) above is   that
the sum   \eqref{recurrentsum}  appearing in the theorem and the measure sums $\sum_{n=1}^N\mu(B(\x,t_n(\x))) $ and $ \sum_{n=1}^N \mu( \hat{R}_n(T,\psi)) $  are all equal up to an additive constant.
   \medskip

  \item[(iii)] The theorem is valid for any self-conformal system on $\R^d$. The price we  seemingly have to pay for this generality is that the radii of the target balls  $B(\x,t_n(\x))$ associated with the modified counting function $ \hat{R}_n(\x,N;T,\psi) $  are dependent on their centres $\x \in K$.  This is   clearly unlike the situation for the ``pure'' counting function  $R_n(\x,N;T,\psi)$   for which we know that Claim~F is false for all self-conformal systems.\medskip
  \item[(iv)]  A simple consequence of Theorem~\ref{quantrec} is the following asymptotic statement that ``fixes'' Claim~F:
 \emph{ Let $(\Phi,K,\mu,T)$ be a self-conformal system on $\R^d$ and let $\psi:\N\to\R_{\geq0}$ be a real positive function such that $\sum_{n=1}^{\infty} \psi(n)$ diverges. Then for $\mu$-almost all $ \x \in  K$}
 \begin{equation}\label{sbccc}
 \lim_{N\to\infty}\frac{\sum_{n=1}^{N}\one_{B(\x,t_n(\x))}
 (T^n\x)}{\sum_{n=1}^{N}\psi(n)} \, = \, 1 \, .
 \end{equation} \\
 Note that in view of the discussion in (ii) above this ``corrected'' statement simply corresponds to Claim~F in which the counting function $R_n(\x,N;T,\psi) $ is replaced by the modified counting function $ \hat{R}_n(\x,N;T,\psi) $    and  $R_n(T,\psi) $ is replaced by $ \hat{R}_n(T,\psi)$.
 \medskip
 \item[(v)] Under various  growth conditions on the function  $\psi$,  Persson \cite{persson2023} has proved a result in a similar  vein to \eqref{sbccc} for a large class of dynamical systems  with exponential decay of correlations on the unit interval. Subsequently, his work (with the various growth conditions) was extended by Rodriguez Sponheimer \cite{Sponheimer2025}  to  more general dynamical systems including Axiom A diffeomorphisms. We stress that Theorem~\ref{quantrec}, which implies \eqref{sbccc},  is free of growth conditions on $\psi$ and provides  an essentially optimal error term.  
Most recently, under a `short return time' assumption and three-fold exponential decay of correlations, Persson and Rodriguez Sponheimer \cite{persson2025}  have essentially removed the growth conditions on $\psi$  that were imposed in their earlier works.
They show that certain non-linear piecewise expanding interval maps of the unit interval, as well   as certain  hyperbolic automorphisms of the two-dimensional torus $\mathbb{T}^2$  satisfy their assumptions.  As a consequence of there main result \cite[Theorem~2]{persson2025}, it then follows that \eqref{sbccc} holds for these dynamical systems.

\end{enumerate}
\end{remark}
\medskip



Even though Theorem~\ref{quantrec}  is in some sense a ``complete'' result, it  fails to  directly deal with the main purpose  of  Claim~F.   Indeed, it remains highly desirable  to  obtain asymptotic information regarding the  behaviour of the ``pure''  counting function \eqref{countdefst} in which the radii of the target balls are independent of their centres.  We reiterate that this is not the case within the framework of Theorem~\ref{quantrec}.   In short, the second main result (for recurrent sets) shows that we are in reasonably  good shape for systems with Gibbs measures equivalent to restricted Hausdorff measures $\cH^{\tau}|_K$.  Here and throughout, we say that Borel measures $\mu$ and $\nu$ on a metric space $(X,d)$ are equivalent if there exists a constant $C\geq1$  such that
$
C^{-1}\nu(E)\leq\mu(E)\leq C\nu(E)  $   for any Borel subset $ E\subseteq X$.  For the proof see \cite[Section~7.2]{jungie}.


\begin{theorem} \label{toprove}
     Let $(\Phi,K,\mu,T)$ be a self-conformal system on $\R^d$ with $\mu$ being a Gibbs measure  equivalent to $\cH^{\tau}|_K$ where $\tau:=\dimH K$. Let $\psi:\N\to\R_{\geq0}$ be a real positive function such that $\psi(n) \to 0 $ as $n \to \infty$. Then for any $\eta>0$ and $\epsilon>0$, we have
     \begin{equation} \label{recurrentcountHM}
        \sum_{n=1}^{N}\one_{B(\x,\psi(n))}(T^n\x)=\sum_{n=1}^N\mu\big(B(\x,\psi(n))\big)+O\left(\Psi_{\eta}(N)^{1/2}(\log\Psi_{\eta}(N))^{\frac{3}{2}+\epsilon}\right)
    \end{equation}
    for $\mu$-almost all $\x\in K$, where
     \begin{equation} \label{recurrentsumHM}
    \Psi_{\eta}(N):=\sum_{n=1}^N\psi(n)^{(1-\eta)\tau}   \, .
    \end{equation}
\end{theorem}

\medskip

\begin{remark}~  \label{newrem} Several comments are in order.

\begin{enumerate}
  \item[(i)]
 Within the setup of self-conformal systems, the notion of  $\mu$ being equivalent to $\cH^{\tau}|_K$  and $\mu$ being $\tau$-Ahlfors regular coincide (see Remark~\ref{HMcoin}).
\\

\item[(ii)] Theorem~\ref{toprove},  could, in principle, be stated with $\eta = \epsilon$ for a cleaner formulation. However, the parameters $\eta$ and $\epsilon$ serve distinct roles. The presence of $\epsilon > 0$ arises necessarily  from the application of a slight generalisation of Lemma~QBC (see \cite[Lemma~7.7]{jungie}) and  cannot be removed.   Notably, this same $\epsilon$ appears in the statements of Theorem~\ref{shrinktarg} and Theorem~\ref{quantrec}.  In contrast, we strongly believe that theorem remains valid with $\eta = 0$, and that the introduction of $\eta>0$ is merely a technical artifact.    We shall soon see that  this  belief is justified  if we are content with asymptotic statements without error term.
  \\

\item[(iii)] In view of (i) it follows that
\begin{equation}  \label{allcomp101}
\sum_{n=1}^N \mu\big(B(\x,\psi(n))\big)   \;  \asymp \;    \sum_{n=1}^N \psi(n)^{\tau} \,  :=  \, \Psi(N)   \, ;
\end{equation}
that is the main term in \eqref{recurrentcountHM} is comparable to  \eqref{recurrentsumHM} with $\eta =0$.  In turn, it follows that due to the presence of $\eta > 0$ in the error term  in \eqref{recurrentcountHM} we can not always conclude that the main term dominates the error term without imposing a condition on the decay rate of $\psi$.  We give a simple example that we hope clearly illustrates the point being made.  Suppose $d=1$ and $\mu$ is one-dimensional Lebesgue measure. For $\alpha > 0$,  consider the function $\psi_\alpha:\N\to\R_{\geq0} $ given by
$$
\psi_{\alpha} (n) \, :=   \,  n^{-\alpha}  \, .
$$
Then, $   \mu\big(B(x,\psi_{\alpha}(n))\big) = 2 n^{-\alpha} $ for any $x \in K$  and so $\sum_{n=1}^\infty  \mu\big(B(x,\psi_{\alpha}(n))\big)$ diverges for any $\alpha \in (0,1]$. Moreover,
\begin{eqnarray*}\label{}
	\lim_{N\to\infty} \; \frac{{\rm Error \ Term \ in \ \eqref{recurrentcountHM} }}{\sum_{n=1}^N\mu\big(B(x,\psi_{\alpha}(n))\big)}  \ = \
\begin{cases}
			0 &\text{if}\ \  \alpha \in (0,1)\\[2ex]
			\infty &\text{if}\ \   \alpha =1 .
		\end{cases}
	\end{eqnarray*}
Thus, Theorem \ref{toprove} does not yield the desired asymptotic statement at the critical exponent $\alpha =1$.    Nevertheless, apart from this flaw, for reasons outlined earlier,  Theorem~\ref{toprove} is on the whole  a more desirable analogue of Theorem~\ref{shrinktarg} than Theorem~\ref{quantrec} for self-conformal systems with Gibbs measures $\mu$ equivalent to  $\cH^{\tau}|_K$.   Under certain additional conditions on the measure (such as $\mu$ being absolutely continuous with respect to  Lebesgue measure) we are able to show that  the flaw  is not present if we are content with asymptotic statements without error terms.
\end{enumerate}
\end{remark}

Note that the theorem shows that  for $\mu$-almost all $\x \in K $,  the asymptotic behaviour  of the counting function $R(\x,N;T,\psi)$  is determined by the behaviour of the measure sum
\begin{equation} \label{recurrentsumHMx}
\sum_{n=1}^N\mu\big(B(\x,\psi(n))\big)  \, ,
\end{equation}
  which, a priori,  is dependent on $\x$.    The point is that  if the measure  $\mu$ is non-uniform,   the measure of the ``target balls''  $B(\x,\psi(n)) $  associated with   $R(\x,N;T,\psi)$  depends on $\x$.
 This is unlike the situation in the shrinking  target framework in which the  measure of the ``target balls''  $B(\y_n,\psi(n)) $  associated with the counting function  $W(\x,N;T,\psi)$ are independent of $\x$.  On a slightly different but related note, we point out that the Gibbs measures associated with the  explicit counterexamples (Examples~\ref{egcantor} \& \ref{egfull})  to Claim~F satisfy the conditions of Theorem~\ref{toprove}.    Thus,    the $\mu$-measure sum \eqref{recurrentsumHMx}   can not in general coincide with the  $\mu$-measure sum involving the sets $R_n(T,\psi)$ associated with the  recurrent $\limsup$ set  $R(T,\psi)$.
However, it is the case  (see \cite[Lemma~7.9]{jungie})  that the sums \eqref{recurrentsumHM} with $\eta =0$,  \eqref{recurrentsumHMx} and $ \sum_{n=1}^N\mu\big(R_n(T,\psi)\big)$ are all comparable\footnote{For the sake of comparison, recall that in the setting of Theorem~\ref{quantrec} the analogous three sums are asymptotically equivalent  (see comment (ii) in Remark~\ref{remarkthm1.5}).}; that is
\begin{equation}  \label{allcomp}
\sum_{n=1}^N \mu\big(B(\x,\psi(n))\big)   \;  \asymp \;   \sum_{n=1}^N \mu\big(R_n(T,\psi)\big) \;  \asymp \;  \Psi(N):= \sum_{n=1}^N \psi(n)^{\tau}  \, .
\end{equation}
Now with Remark~\ref{remarkthm1.5}\,(ii) in mind,  it follows that if $\mu$ is  $\tau$-Ahlfors regular  then for all $n \in \N$
\begin{equation}  \label{allcomp2}
 \hat{R}_n\big(C^{-1} \psi(n)^{\tau}\big)
 \; \subseteq   \;  R_n(T,\psi) \;   \subseteq  \;  \hat{R}_n\big(C \psi(n)^{\tau}\big)    \, ,
  \end{equation}
  where $C \ge 1 $ is the ``Ahlfors regular'' constant appearing in \eqref{tar}.  Then on making use of \eqref{allcomp} and
  \eqref{allcomp2},  it is easily verified that  Theorem~\ref{quantrec} implies  the   following  zero-one measure criterion  which validates \eqref{divcondsv1} whenever $\mu$ is equivalent to $\cH^{\tau}|_K$.  Indeed, it  coincides with the main result of  Baker $\&$ Farmer \cite{baker2021} discussed within the context of~\eqref{divcondsv1}.


   \begin{corollary} \label{toprovecor}
     Let $(\Phi,K,\mu,T)$ be a self-conformal system on $\R^d$ with $\mu$ being a Gibbs measure equivalent to $\cH^{\tau}|_K$ where $\tau:=\dimH K$. Let $\psi:\N\to\R_{\geq0}$ be a real positive function.  Then \begin{eqnarray}\label{0-1law}
		 \mu\left(R(T,\psi)\right)=
		\begin{cases}
			0 &\text{if}\ \  \sum_{n=1}^{\infty}\psi(n)^{\tau} <\infty,\\[2ex]
			1 &\text{if}\ \ \sum_{n=1}^{\infty} \psi(n)^{\tau}  =\infty.
		\end{cases}
	\end{eqnarray}
\end{corollary}

  \noindent As we have seen, the corollary follows in a fairly straightforward manner from Theorem~\ref{quantrec}. However, it would follow in an entirely trivial way if we were able to eliminate the dependence on $\eta > 0$ in the statement of Theorem~\ref{toprove}. At present, this is not possible, and so we must invoke Theorem~\ref{quantrec} as stated.

We now point out that in the case $\mu$ is equivalent to $\cH^{\tau}|_K$, beyond  implying  the above zero-one measure criterion,   Theorem~\ref{quantrec} can  also  be utilized to explicitly obtain information regarding the behaviour of the counting function \eqref{countdefst}.  In order to state precisely what exactly can be obtained,  we need to introduce the following notion of upper and lower densities.  Let $\psi:\N\to\R_{\geq0}$ be a real positive function.   Then, for each $\tau>0$, each probability measure $\mu$ on $\R^d$ and each $\x\in\R^d$,  we define the $\tau$-lower and $\tau$-upper densities of $\mu$ at $\x$ associated with $\psi$ by
\[
\Theta^{\tau}_*(\mu, \psi,\x):=\liminf_{n\to \infty}\frac{\mu(B(\x,\psi(n)))}{\psi(n)^{\tau}} \, ,  \quad \Theta^{*\tau}(\mu, \psi,\x):=\limsup_{n\to\infty}\frac{\mu(B(\x,\psi(n)))}{\psi(n)^{\tau}}.
\]
With this in mind, the following can be deduced directly from Theorem~\ref{quantrec} - see \cite[Section~7.3]{jungie} for the details.

\begin{theorem} \label{quantcount}
   Let $(\Phi,K,\mu,T)$ be a self-conformal system on $\R^d$ with $\mu$ being a Gibbs measure equivalent to $\cH^{\tau}|_K$ where $\tau:=\dimH K$. Let $\psi:\N\to\R_{\geq0}$ be a real positive function such that $\psi(n) \to 0 $ as $n \to \infty$ and assume that $\sum_{n=1}^{\infty}\psi(n)^{\tau}
    $ diverges.  Then, for $\mu$-almost  all $ \x \in  K$
    \begin{eqnarray*}\label{hauscount}
    \begin{split}
\frac{\Theta^{\tau}_*(\mu,\psi,\x)}{\Theta^{*\tau}(\mu,\psi,\x)}
&\leq\liminf_{N\to\infty}
\frac{\sum_{n=1}^N\one_{B(\x,\psi(n))}(T^n\x)}{\sum_{n=1}^N\mu(B(\x,\psi(n)))}\\
&\leq\limsup_{N\to\infty}
\frac{\sum_{n=1}^N\one_{B(\x,\psi(n))}(T^n\x)}{\sum_{n=1}^N\mu(B(\x,\psi(n)))}
\leq\frac{\Theta^{*\tau}(\mu,\psi,\x)}{\Theta^{\tau}_*(\mu,\psi,\x)}  \, .
    \end{split}
    \end{eqnarray*}
\end{theorem}

 \medskip

 \noindent Clearly, this result is  weaker  than Theorem~\ref{toprove} whenever the main term in \eqref{recurrentcountHM} dominates. In such cases, Theorem~\ref{toprove} yields an asymptotic statement with a quantified error term, whereas Theorem~\ref{quantcount} provides, at best, an unquantified asymptotic statement. However, when the error term in \eqref{recurrentcountHM} dominates, Theorem~\ref{toprove} becomes ineffective, while Theorem~\ref{quantcount} may still yield meaningful information—and is potentially stronger. In particular, when the lower and upper densities of $\mu$ associated with $\psi$ coincide, we can derive the following asymptotic statement directly from Theorem~\ref{quantcount}.


\begin{corollary} \label{quantcountcor1}
   Let $(\Phi,K,\mu,T)$ be a self-conformal system on $\R^d$ with $\mu$ being a Gibbs measure equivalent to $\cH^{\tau}|_K$ where $\tau:=\dimH K$. Let $\psi:\N\to\R_{\geq0}$ be a real positive function such that $\psi(n) \to 0 $ as $n \to \infty$ and assume that $\sum_{n=1}^{\infty}\psi(n)^{\tau}
    $ diverges and that
     $ \Theta^{\tau}_*(\mu, \psi,\x) =  \Theta^{*\tau}(\mu, \psi,\x) $
    for $\mu$-almost  all $ \x \in  K$.    Then, for $\mu$-almost  all $ \x \in  K$
    \begin{eqnarray*}
    \begin{split}
\lim_{N\to\infty}
\frac{\sum_{n=1}^N\one_{B(\x,\psi(n))}(T^n\x)}{\sum_{n=1}^N\mu(B(\x,\psi(n)))} =1
    \end{split}
    \end{eqnarray*}
\end{corollary}

\medskip

We now consider the special case in which the Gibbs measure is absolutely continuous with respect to $d$-dimensional Lebesgue measure $\cL^d$.  For convenience,  let   $c_d:=\cL^d\big(B(0,1)\big)$ and suppose that $\mu$ is a Gibbs measure equivalent to $\cL^{d}|_K$ with density function $h$.
Then,   the  Lebesgue density theorem
implies  that for $\mu$-almost all $\x\in K$
\smallskip
\begin{equation}  \label{ldt}
    \mu\big(B(\x,\psi(n))\big)= \big(h(\x)+\epsilon_n(\x)\big)\cdot c_d  \, \psi(n)^d  \, ,
\end{equation}

\noindent where $\epsilon_n(\x)\to 0$ as $n\to\infty$.   The upshot of this is the following statement for  absolutely continuous  measures.  The first part is  a rewording  of Theorem~\ref{toprove} while the second part is a  rewording of  Corollary~\ref{quantcountcor1}.


\begin{corollary}
\label{quantcountabcont}
     Let $(\Phi,K,\mu,T)$ be a self-conformal system on $\R^d$ and suppose that $\dimH K=d$. Let  $\mu$  be a  Gibbs measure equivalent to $\cL^{d}|_K$ with density function $h$.  Let $\psi:\N\to\R_{\geq0}$ be a real positive function such that $\psi(n) \to 0 $ as $n \to \infty$. Then the following are true.
\begin{itemize}
  \item[(i)] For any $\eta>0$ and $\epsilon>0$, we have
  \begin{eqnarray} \label{lebcount}
        \sum_{n=1}^{N}\one_{B(\x,\psi(n))}(T^n\x) \;
        &=&  \;   c_d h(\x)\Psi(N)   +   c_d\sum_{n=1}^N\epsilon_n(\x)\psi(n)^d  \nonumber  \\[1ex]  &~&  \hspace*{10ex} +\ O\left(\Psi_{\eta}(N)^{1/2}(\log\Psi_{\eta}(N))^{\frac{3}{2}+\epsilon}\right).
    \end{eqnarray}
    for $\mu$-almost all $\x\in K$, where   $\Psi_{\eta}(N) :=  \sum_{n=1}^N\psi(n)^{(1-\eta)d}  \, $,   $\Psi(N) := \Psi_{0}(N)  \, $  and  $\epsilon_n(\x)\to 0$ as $n\to\infty$ satisfies~\eqref{ldt}.\\

  \item[(ii)]  If  \ $\sum_{n=1}^\infty \psi(n)^{d}$ diverges,  then
   \begin{equation}  \label{ghjk}
\lim_{n \to \infty} \frac{\sum_{n=1}^{N}\one_{B(\x,\psi(n))}(T^n\x) }{c_d\Psi(N) }  \, =  \, h(\x)  \quad \text{for $\mu$-almost all $\x\in K$.  }    \end{equation}
\end{itemize}

\end{corollary}


Note that  in general we do not have any information regarding the rate at which  $\epsilon_n(\x)\to0$, so it is not possible to compare the size of the second and third terms appearing on the right hand side of
\eqref{lebcount}.  However, if $\mu= \cL^{d}$ then $\epsilon_n(\x)= 0$ for all $n \in \N $ and $\x \in K$ and so the second term is zero.
With this in mind, it follows that Corollary~\ref{quantcountabcont} is in line with the main result established in \cite{levesley2024} for piecewise linear maps of $[0,1]^d$.   Furthermore, with the fact that self-conformal  systems are exponentially mixing at our disposal (the key result established in \cite{jungie}),  the   asymptotic statement  \eqref{ghjk}
can be  directly derived from the recent impressive  work of He \cite{he2024}. In short, He obtains \eqref{ghjk}  for a class of  measure-preserving systems  for which $\mu$ is exponentially mixing and absolutely continuous with respect to  Lebesgue measure.

\bigskip

We bring this section to an end with a brief discussion concerning  the  recurrent problem  beyond self-conformal systems, or rather beyond the structure inherited by such systems.    In view of Theorem~\ref{shrinktarg}, we know that  exponential mixing underpins the asymptotic behaviour of the counting function within the setup of the shrinking target problem. Currently, we see no obvious counterexample that  shows that this is not enough within the recurrent framework.  Adding a safety net, by  restricting to Hausdorff  measures, it remains plausible that the following  ``strengthening''  of Theorem~\ref{toprove} is true.   In short it would suggest  that the  key aspect of the system under consideration is that it is exponentially mixing and nothing else.

\noindent\textbf{Claim T.  \ }
 \emph{Let $(X,\mathcal{B},\mu,T)$ be a measure-preserving dynamical system in $\R^d$  with $\mu$ being a  $\tau$-Ahlfors regular measure where $\tau:=\dimH X$.   Suppose that $T$ is exponentially mixing with respect to $\mu$. Let $\psi:\N\to\R_{\geq0}$ be a real positive function such that $\psi(n) \to 0 $ as $n \to \infty$.  Then, for any given $\varepsilon>0$, we have that
 \begin{equation*}
        \sum_{n=1}^N\one_{B(\x,\psi(n))}(T^n\x)=\Psi(N,\x)+O\left(\Psi(N,\x)^{1/2}\log^{\frac{3}{2}+\epsilon}(\Psi(N,x))\right)
    \end{equation*}
    for $\mu$-almost all $\x\in X$, where
$\Psi(N,\x):=\displaystyle{\sum_{n=1}^N}\mu\big(B(\x,\psi(n))\big)  $\,.  }

\bigskip

\noindent Several comments are in order.
\begin{enumerate}
\item[(i)]  Recall Remark~\ref{newrem}\,(i), namely that within the setup of self-conformal systems, the notion of  $\mu$ being equivalent to $\cH^{\tau}|_X$ and $\mu$ being $\tau$-Ahlfors regular coincide.
\medskip

\item[(ii)] Clearly, under the assumption that  $\mu$ is a  $\tau$-Ahlfors regular measure as in Claim~T,  we can replace the quantity  $\Psi(N,\x)$ by $\sum_{n=1}^N\psi(n)^{\tau}  $  in the error term  and thus making it independent of $ \x \in X$.    The reason that we have not done this is that there is a possibility that the conclusion of the claim is true without the Ahlfors regular assumption and in such generality the error may depend on $\x \in X$; that is to say that  $\Psi(N,\x)$ may not be comparable to a sum that is independent of $\x$.
\medskip

\item[(iii)] With the previous comment in mind, it is worth pointing out that \eqref{allcomp}  is in fact true under the hypothesis of Claim~T  (see \cite[Lemma~2.5]{hussain2022}).
Indeed, it is easily checked that all that is essentially required to  establish \eqref{allcomp}  is that $\mu$ is $\tau$-Ahlfors regular and that $\mu$ is exponentially mixing.
\end{enumerate}


\noindent Even if Claim~T turns out to be false, it does not rule out the following strengthening of Corollary~\ref{toprovecor} which is of independent interest.

\noindent\textbf{Claim 0-\!1.}  \emph{Let $(X,\mathcal{B},\mu,T)$ be a measure-preserving dynamical system in $\R^d$  with $\mu$ being a  $\tau$-Ahlfors regular measure where $\tau:=\dimH X$.  Suppose that $T$ is exponentially mixing with respect to $\mu$. Let $\psi:\N\to\R_{\geq0}$ be a real positive function.  Then \begin{eqnarray}\label{claim0-1law}
		 \mu\left(R(T,\psi)\right)=
		\begin{cases}
			0 &\text{if}\ \  \sum_{n=1}^{\infty}\psi(n)^{\tau} <\infty,\\[2ex]
			1 &\text{if}\ \ \sum_{n=1}^{\infty} \psi(n)^{\tau}  =\infty.
		\end{cases}
	\end{eqnarray}}

\medskip

\noindent As already mentioned,  currently we see no obvious counterexample that  shows that Claim~T is false, let alone a counterexample to Claim 0-\!1.

\begin{remark}  \label{nonono}
As mentioned in the discussion leading up to Claim~T, the actual statement of the claim is erring on the side of caution. Indeed, we see no obvious counterexample to either Claim~T or Claim~0-1 even if we remove the assumption that the measure $\mu$ is Ahlfors regular.  Obviously, without the latter assumption,  in Claim~0-1  the conclusion  \eqref{claim0-1law}  would read:
\begin{eqnarray*}
		 \mu\left(R(T,\psi)\right)=
		\begin{cases}
			0 &\text{if}\ \  \sum_{n=1}^{\infty}\mu\big(B(\x,\psi(n))\big) <\infty   \quad \text{for  } \text{$\mu$-almost all \ } \x \in X \, ,  
            \\[2ex]
			1 &\text{if}\ \ \sum_{n=1}^{\infty} \mu\big(B(\x,\psi(n))\big)  =\infty   \quad \text{for  } \text{$\mu$-almost all \ } \x \in X \,   .
		\end{cases}
	\end{eqnarray*}
It is worth pointing out that   a relatively painless  calculation shows that within the context of Example~ABB, we have that
$$
\textstyle{\sum_{n=1}^{\infty}\mu\big(B(\x,\psi_{\alpha}(n))\big) \, <  \, \infty}
$$
for $\mu$-almost all $\x \in K$ (see  \cite[Appendix~C]{jungie}  for the details).
Thus,  Example~ABB is not a counterexample to the bolder statement in which the Ahlfors regular assumption is dropped.  Finally,  at the very basic level,  as far as we are aware, it is not known whether or not $\mu(R(T,\psi))$ satisfies a zero-one law; i.e.  $\mu(R(T,\psi))= 0 $ or $1$.
\end{remark}

%



\bibliographystyle{abbrv}

\bibliography{bibliography}

\begin{thebibliography}{10}

\bibitem{allen2025}
D.~Allen, S.~Baker, and B.~B{\'a}r{\'a}ny.
\newblock Recurrence rates for shifts of finite type.
\newblock {\em Adv. Math.}, 460:110039, 2025.

\bibitem{allen2021}
D.~Allen and B.~B{\'a}r{\'a}ny.
\newblock On the {H}ausdorff measure of shrinking target sets on self-conformal
  sets.
\newblock {\em Mathematika}, 67(4):807--839, 2021.

\bibitem{aspenberg2019}
M.~Aspenberg and T.~Persson.
\newblock Shrinking targets in parametrised families.
\newblock {\em Math. Proc. Cambridge Philos. Soc.}, 166(2):265--295, 2019.

\bibitem{baker2021}
S.~Baker and M.~Farmer.
\newblock Quantitative recurrence properties for self-conformal sets.
\newblock {\em Proc. Amer. Math. Soc.}, 149(3):1127--1138, 2021.

\bibitem{baker2024}
S.~Baker and H.~Koivusalo.
\newblock Quantitative recurrence and the shrinking target problem for
  overlapping iterated function systems.
\newblock {\em Adv. Math.}, 442:109538, 2024.

\bibitem{barany2018}
B.~B{\'a}r{\'a}ny and M.~Rams.
\newblock Shrinking targets on {B}edford--{M}c{M}ullen carpets.
\newblock {\em Proc. Lond. Math. Soc.}, 117(5):951--995, 2018.

\bibitem{RothVol}
V.~Beresnevich, V.~Bernik, M.~Dodson, S.~Velani, W.~Chen, T.~Gowers,
  H.~Halberstam, W.~Schmidt, and R.~Vaughan.
\newblock {\em Classical metric Diophantine approximation revisited}, pages
  38--61.
\newblock Cambridge University Press, United States, Feb. 2009.

\bibitem{BHV2024}
V.~Beresnevich, M.~Hauke, and S.~Velani.
\newblock Borel--{C}antelli, zero-one laws and inhomogeneous
  {D}uffin--{S}chaeffer.
\newblock {\em arXiv preprint arXiv:2406.19198}, 2024.

\bibitem{BRV2016}
V.~Beresnevich, F.~Ram\'{\i}rez, and S.~Velani.
\newblock Metric {D}iophantine approximation: aspects of recent work.
\newblock In {\em Dynamics and analytic number theory}, volume 437 of {\em
  London Math. Soc. Lecture Note Ser.}, pages 1--95. Cambridge Univ. Press,
  Cambridge, 2016.

\bibitem{BeresnevichVelani7}
V.~Beresnevich and S.~Velani.
\newblock The divergence {B}orel-{C}antelli lemma revisited.
\newblock {\em J. Math. Anal. Appl.}, 519(1):Paper No. 126750, 21, 2023.

\bibitem{bill}
P.~Billingsley.
\newblock {\em Probability and Measure}.
\newblock Wiley Series in Probability and Statistics. Wiley, 1995.

\bibitem{borel}
A.~Borel.
\newblock {\em Linear {A}lgebraic {G}roups, second enlarged edition}, volume
  126 of {\em Graduate Texts in Mathematics}.
\newblock Springer, 1991.

\bibitem{boshernitzan1993}
M.~D. Boshernitzan.
\newblock Quantitative recurrence results.
\newblock {\em Invent. Math.}, 113(1):617--631, 1993.

\bibitem{bugeaud2014}
Y.~Bugeaud and B.-W. Wang.
\newblock Distribution of full cylinders and the {D}iophantine properties of
  the orbits in $\beta$-expansions.
\newblock {\em J. Fractal Geom.}, 1(2):221--241, 2014.

\bibitem{chang2019}
Y.~Chang, M.~Wu, and W.~Wu.
\newblock Quantitative recurrence properties and homogeneous self-similar sets.
\newblock {\em Proc. Amer. Math. Soc.}, 147(4):1453--1465, 2019.

\bibitem{fan1999}
A.~H. Fan and K.-S. Lau.
\newblock Iterated function system and {R}uelle operator.
\newblock {\em J. Math. Anal. Appl.}, 231(2):319--344, 1999.

\bibitem{fang2020}
L.~Fang, M.~Wu, and B.~Li.
\newblock Approximation orders of real numbers by $\beta$-expansions.
\newblock {\em Math. Z.}, 296(1):13--40, 2020.

\bibitem{fernandez2012}
J.~Fern{\'a}ndez, M.~Meli{\'a}n, and D.~Pestana.
\newblock Expanding maps, shrinking targets and hitting times.
\newblock {\em Nonlinearity}, 25(9):2443, 2012.

\bibitem{galatolo2015}
S.~Galatolo, J.~Rousseau, and B.~Saussol.
\newblock Skew products, quantitative recurrence, shrinking targets and decay
  of correlations.
\newblock {\em Ergodic Theory Dynam. Systems}, 35(6):1814--1845, 2015.

\bibitem{ghosh2024}
A.~Ghosh and D.~Nandi.
\newblock {D}iophantine approximation, large intersections and geodesics in
  negative curvature.
\newblock {\em Proc. Lond. Math. Soc.}, 128(2):e12581, 2024.

\bibitem{Harman1998}
G.~Harman.
\newblock {\em Metric {N}umber {T}heory}.
\newblock London Mathematical Society Monographs. New Series, 18. The Clarendon
  Press, Oxford University Press, New York, 1998.

\bibitem{he2024}
Y.~He.
\newblock Quantitative recurrence properties and strong dynamical
  $\mathrm{B}$orel-$\mathrm{C}$antelli lemma for dynamical systems with
  exponential decay of correlations.
\newblock {\em arXiv preprint arXiv:2410.10211}, 2024.

\bibitem{he2022}
Y.~He and L.~Liao.
\newblock Quantitative recurrence properties for piecewise expanding maps on
  $[0, 1]^{d}$.
\newblock {\em Ann. Sc. Norm. Super. Pisa Cl. Sci.}, pages 1--40, 2024.

\bibitem{hill1995}
R.~Hill and S.~L. Velani.
\newblock The ergodic theory of shrinking targets.
\newblock {\em Invent. Math.}, 119(1):175--198, 1995.

\bibitem{MR1471868}
R.~Hill and S.~L. Velani.
\newblock Metric {D}iophantine approximation in {J}ulia sets of expanding
  rational maps.
\newblock {\em Inst. Hautes \'Etudes Sci. Publ. Math.}, 85:193--216, 1997.

\bibitem{hu2024}
Z.~N. Hu and T.~Persson.
\newblock Hausdorff dimension of recurrence sets.
\newblock {\em Nonlinearity}, 37(5), 2024.

\bibitem{jungie}
J.~Huang, B.~Li, and S.~Velani.
\newblock Exponential mixing for gibbs measures on self-conformal sets and
  applications.
\newblock {\em arXiv preprint arXiv:2504.00632}, 2025.

\bibitem{hussain2022}
M.~Hussain, B.~Li, D.~Simmons, and B.~Wang.
\newblock Dynamical {B}orel--{C}antelli lemma for recurrence theory.
\newblock {\em Ergodic Theory Dynam. Systems}, 42(6):1994--2008, 2022.

\bibitem{khintchine1924}
A.~Khintchine.
\newblock {E}inige {S}{\"a}tze {\"u}ber {K}ettenbr{\"u}che, mit {A}nwendungen
  auf die {T}heorie der diophantischen {A}pproximationen.
\newblock {\em Math. Ann.}, 92(1):115--125, 1924.

\bibitem{kim2007}
D.~H. Kim.
\newblock The shrinking target property of irrational rotations.
\newblock {\em Nonlinearity}, 20(7):1637, 2007.

\bibitem{Kirsebom2023}
M.~Kirsebom, P.~Kunde, and T.~Persson.
\newblock On shrinking targets and self-returning points.
\newblock {\em Ann. Sc. Norm. Super. Pisa Cl. Sci. (5)}, 24(3):1499--1535,
  2023.

\bibitem{kleinbock2023}
D.~Kleinbock and J.~Zheng.
\newblock Dynamical {B}orel--{C}antelli lemma for recurrence under {L}ipschitz
  twists.
\newblock {\em Nonlinearity}, 36(2):1434, 2023.

\bibitem{Kurzweil2}
J.~Kurzweil.
\newblock On the metric theory of inhomogeneous diophantine approximations.
\newblock {\em Studia Math.}, 15:84--112, 1955.

\bibitem{levesley2024}
J.~Levesley, B.~Li, D.~Simmons, and S.~Velani.
\newblock Shrinking targets versus recurrence: the quantitative theory.
\newblock {\em Mathematika}, 71(4):Paper No. e70039, 16, 2025.

\bibitem{LI2023108994}
B.~Li, L.~Liao, S.~Velani, and E.~Zorin.
\newblock The shrinking target problem for matrix transformations of tori:
  Revisiting the standard problem.
\newblock {\em Adv. Math.}, 421:108994, 2023.

\bibitem{li2014}
B.~Li, B.-W. Wang, J.~Wu, and J.~Xu.
\newblock The shrinking target problem in the dynamical system of continued
  fractions.
\newblock {\em Proc. Lond. Math. Soc.}, 108(1):159--186, 2014.

\bibitem{persson2019}
T.~Persson.
\newblock Inhomogeneous potentials, {H}ausdorff dimension and shrinking
  targets.
\newblock {\em Ann. H. Lebesgue}, 2:1--37, 2019.

\bibitem{persson2023}
T.~Persson.
\newblock A strong $\mathrm{B}$orel--$\mathrm{C}$antelli lemma for recurrence.
\newblock {\em Studia Math.}, 268(1):75--89, 2023.

\bibitem{persson2025}
T.~Persson and A.~Rodriguez~Sponheimer.
\newblock Strong {B}orel--{C}antelli {L}emmas for {R}ecurrence.
\newblock {\em Accepted for publication in Studia Math.}, 2025.

\bibitem{Port}
S.~C. Port.
\newblock {\em Theoretical probability for applications}.
\newblock Wiley Series in Probability and Mathematical Statistics: Probability
  and Mathematical Statistics. John Wiley \& Sons, Inc., New York, 1994.
\newblock A Wiley-Interscience Publication.

\bibitem{Sponheimer2025}
A.~Rodriguez~Sponheimer.
\newblock A recurrence-type strong {B}orel-{C}antelli lemma for {A}xiom {A}
  diffeomorphisms.
\newblock {\em Ergodic Theory Dynam. Systems}, 45(3):936--955, 2025.

\bibitem{seuret2015}
S.~Seuret and B.-W. Wang.
\newblock Quantitative recurrence properties in conformal iterated function
  systems.
\newblock {\em Adv. Math.}, 280:472--505, 2015.

\bibitem{Sprindzuk}
V.~G. Sprind{\v z}uk.
\newblock {\em Metricheskaya teoriya diofantovykh priblizheni\u\i~ ({M}etric
  theory of {D}iophantine approximations}.
\newblock Izdat. ``Nauka'', Moscow, 1977.

\bibitem{tan2011}
B.~Tan and B.-W. Wang.
\newblock Quantitative recurrence properties for beta-dynamical system.
\newblock {\em Adv. Math.}, 228(4):2071--2097, 2011.

\bibitem{tseng2008}
J.~Tseng.
\newblock On circle rotations and the shrinking target properties.
\newblock {\em Discrete Contin. Dyn. Syst.}, 20(4):1111--1122, 2008.

\bibitem{walters2000}
P.~Walters.
\newblock {\em An introduction to ergodic theory}, volume~79.
\newblock Springer Science \& Business Media, 2000.

\bibitem{wu2025}
Y.-L. Wu and N.~Yuan.
\newblock Quantitative recurrence problem on some {B}edford-{M}c{M}ullen
  carpets.
\newblock {\em J. Math. Anal. Appl.}, 543(2):Paper No. 128938, 15, 2025.

\end{thebibliography}

\end{document}